\documentclass[11pt,a4paper,fleqn]{article}
\RequirePackage{amsthm,amsmath,natbib}
\RequirePackage{amssymb,amstext}
\RequirePackage{hypernat}
\usepackage[pdfpagemode=UseOutlines ,plainpages=false,hypertexnames=false ,pdfpagelabels ,hyperindex=true,colorlinks=true]{hyperref}
	\makeatletter%
	\Hy@breaklinkstrue%
	\makeatother%
\usepackage{color}
\definecolor{darkred}{rgb}{0.6,0.0,0.1}
\definecolor{darkgreen}{rgb}{0,0.5,0}
\definecolor{darkblue}{rgb}{0,0,0.5}
\hypersetup{colorlinks ,linkcolor=darkblue ,filecolor=darkgreen
,urlcolor=darkblue ,citecolor=black}

\usepackage{		url}
\usepackage{verbatim}
\usepackage{ bm}

\renewcommand{\cite}{\citet}
\usepackage{jan-abrev-package}

\definecolor{dgreen}{rgb}{0,0.5,0}
\definecolor{dblue}{rgb}{0,0,0.9}
\definecolor{dred}{rgb}{0.6,0.0,0.1}
\definecolor{dgold}{rgb}{0.5,0.3,0.0}
\definecolor{dvio}{rgb}{0.6,0.3,0.5}
\definecolor{gray}{rgb}{0.5,0.5,0.5}

\newtheoremstyle{mysc}
  {3pt}
  {3pt}
  {\it}
  {}
  {\color{darkred}\sc}
  {.}
  {.5em}
  {}

\newtheoremstyle{myex}
  {10pt}
  {10pt}
  {\rm}
  {}
  {\color{darkred}\sc}
  {.}
  {.5em}
  {}

\theoremstyle{mysc}\newtheorem{prop}{Proposition}[section]
\theoremstyle{mysc}\newtheorem{assumption}{Assumption}[section]
\theoremstyle{mysc}\newtheorem{coro}[prop]{Corollary}
\theoremstyle{mysc}\newtheorem{theo}[prop]{Theorem}
\theoremstyle{mysc}
\theoremstyle{mysc}\newtheorem{lem}[prop]{Lemma}
\theoremstyle{myex}\newtheorem{rem}{Remark}[section]
\theoremstyle{myex}
\theoremstyle{myex}

\numberwithin{equation}{section}

\makeatletter%
\def\@fnsymbol#1{\ensuremath{\ifcase#1\or * \or \star \or 1 \or 2\or 3\or  , \or
g\or h\or i\else\@ctrerr\fi}}%
\makeatother%

\author{{\sc Fabienne Comte}\thanks{Universit\'e Paris Descartes, Laboratoire MAP5, UMR CNRS 8145, 45, rue des Saints-P\`eres, F-75270 Paris cedex 06, France, e-mail:  \url{fabienne.comte@parisdescartes.fr}} \and {\sc Jan Johannes}\thanks{Universit\"at Heidelberg, Institut f\"ur Angewandte Mathematik, Im Neuenheimer Feld, 294, D-69120 Heidelberg, Germany, e-mail: \url{johannes@statlab.uni-heidelberg.de}}}

\title{{\bf Adaptive estimation in circular functional linear models.}}

\begin{document}
\date{\today}
\maketitle

\begin{abstract} We consider the problem of estimating the slope parameter in  circular functional linear regression, where scalar responses
$Y_1,\dotsc,Y_n$ are modeled in dependence of $1$-periodic, second order stationary random functions
$X_1,\dotsc,X_n$.  We consider an orthogonal series estimator of the slope function $\beta$, by replacing the first  $m$ theoretical coefficients of its development in the trigonometric basis by adequate estimators. We propose a model selection procedure for $m$ in a set of admissible values, by defining a contrast function minimized by our estimator and a theoretical penalty function; this first step assumes the degree of ill posedness to be known. Then we generalize the procedure to a random set of admissible $m$'s and a random penalty function. The resulting estimator is completely data driven and reaches automatically what is known to be the optimal minimax rate of convergence, in term of a general weighted $L^2$-risk. This means that we provide adaptive estimators of both $\beta$ and its derivatives.\end{abstract}

\begin{tabbing}
\noindent \emph{Keywords:} \=Orthogonal series estimation; model selection; derivatives estimation;\\
\> mean squared error of prediction; minimax theory.\\[.2ex]
\noindent\emph{AMS 2000 subject classifications:} Primary 62G05; secondary 62J05, 62G08.
\end{tabbing}

\section{Introduction}\label{sec:intro}
Functional linear models have become very important in a diverse range of disciplines,
including medicine, linguistics, chemometrics as well as econometrics (see for instance 
\cite{RamsaySilverman2005} and \cite{FerratyVieu2006}, for several case studies, or more specific, \cite{ForniReichlin1998} and \cite{PredaSaporta2005} for  applications in economics).   Roughly speaking, in all
these applications the dependence 
of a response variable $Y$  on the variation of an explanatory random function $X$ is modeled by 
\begin{equation}\label{intro:e1}Y=\int_{0}^1\sol(t)X(t)dt +\sigma\epsilon,\quad\sigma>0,
\end{equation}
for some error term $\epsilon$. One objective is then to estimate nonparametrically the
slope function $\sol$  based on an independent and identically distributed (i.i.d.) sample of 
$(Y,X)$.  

In this paper we suppose  that the random function $X$ is taking  its  values in $L^2[0,1]$, which
is endowed with the usual inner product $\skalar$ and induced norm $\norm$, and that $X$ has a finite second moment, i.e., $\Ex \normV{X}^2<\infty$.  In order to simplify notations we assume that the mean function of $X$ is zero. Moreover, the random function $X$ and the error term $\epsilon$ are uncorrelated, where  $\epsilon$ is assumed to have  mean zero and variance one. This situation
has been considered, for example, in \cite{CardotFerratySarda2003}, \cite{MullerStadtmuller2005} or most recently \cite{AOStoappear}. Then multiplying both sides in (\ref{intro:e1}) by $X(s)$ and taking the expectation leads to 
\begin{equation}\label{method:e1}g(s):=\Ex[YX(s)]=
\int_{0}^1\beta(t)\cov(X(t),X(s))dt=:[\op \sol](s),\quad s\in[0,1],
\end{equation}
where $g$ belongs to $L^2[0,1]$ and $\op$ denotes the covariance operator associated to the random function $X$.  We shall assume that there exists a unique solution $\sol\in L^2[0,1]$ of equation (\ref{method:e1}). Estimation of $\sol$ is thus linked with the
inversion of the covariance operator $\op$  and, known to be an ill-posed inverse problem (for a detailed discussion in the context of inverse problems see chapter 2.1 in \cite{EHN00}, while in the special case of a functional linear model we refer to \cite{CardotFerratySarda2003}). 

In this paper we consider a circular functional linear model (defined below), where  the associated covariance operator $\op$ admits a spectral decomposition $\{\ev_j,\ef_j,j\geq1\}$ given by the trigonometric basis $\{\ef_j\}$ as eigenfunctions and a strictly positive, possibly not ordered, zero-sequence $\ev:=(\ev_j)_{j\geq 1}$ of corresponding eigenvalues.  Then the normal equation can be rewritten as follows 
\begin{equation}\label{intro:e2}
\sol=\sum_{j=1}^\infty \frac{[g]_j}{\ev_j}\cdot\ef_j\quad\text{ with } [g]_j:=\skalarV{{g},\ef_j},\; j\geq1.
\end{equation} 
For estimation purpose, we replace the unknown quantities $g_j$ and $\ev_j$ in equation \eqref{intro:e2} by their empirical counterparts. That is, if  $(Y_1,X_1),\dotsc,(Y_n,X_n)$ denotes an i.i.d. sample of $(Y,X)$, then for each $j\geqslant 1$, we consider the unbiased estimator
\begin{equation*}
[\widehat{g}]_j:=\frac{1}{n}\sum_{i=1}^n Y_{i}\,[X_i]_j,\quad\mbox{and}\quad \hev_{j}:=\frac{1}{n}\sum_{i=1}^n[X_i]_j^2\quad\mbox{ with } [X_i]_j:=\skalarV{X_{i},\ef_{j}}\end{equation*}
for $[g]_j$ and $\ev_j$ respectively. The orthogonal series estimator $\hsol_m$ of $\sol$  is then defined  by
\begin{equation}\label{intro:def:est:reg}
\hsol_m:=\sum_{j=1}^m \frac{\widehat{g}_{j}}{\hev_{j}} \cdot\1\{\hev_{j}\geqslant 1/n\}\cdot \ef_{j}.\end{equation} 
Note that we introduce an additional threshold $1/n$ on each estimated eigenvalue $\hev_{j}$, since it could be arbitrarily close to zero  even in case that the true eigenvalue $\ev_j$ is sufficiently far away from zero. Moreover, the orthogonal series estimator keeps only $m$ coefficients; this is an alternative to the popular  Tikhonov regularization (c.f. \cite{HallHorowitz2007}), where in (\ref{intro:e2})  the factor $1/\ev_j$ is replaced by $\ev_j/(\alpha+\ev_j^2)$.  Thresholding in the Fourier domain  has been used, for example, in a deconvolution problem in  \cite{MairRuymgaart96} or \cite{Neumann1997} and coincides with an approach called spectral cut-off in the numerical analysis literature  (c.f. \cite{Tautenhahn96}).  
 
In this paper we  shall measure the performance of an estimator $\hsol$ of $\beta$  by the $\cF_{\hw}$-risk, that is $\Ex\|\hsol-\sol\|_{\hw}^2$, where for  some strictly positive sequence of weights  $\hw:=(\hw_j)_{j\geq 1}$
\begin{equation*}
\normV{f}_\hw^2 := \sum_{j=1}^\infty \hw_j |\skalarV{f,\ef_j}|^2\qquad\mbox{ for all $f \in L^2[0,1].$ }
\end{equation*}
This general framework allows us with appropriate choices of the weight sequence $\hw$ to cover the estimation not only of the slope parameter itself (c.f. \cite{HallHorowitz2007}) but also of its derivatives  as well as the optimal estimation with respect to the mean squared prediction error (c.f. \cite{CardotFerratySarda2003} or \cite{CrambesKneipSarda2007}). For a more detailed discussion, we refer to \cite{CardotJohannes2007}. It is well-known that the obtainable accuracy of any estimator in terms of the $\cF_\hw$-risk  is essentially determined by the regularity conditions imposed on both the slope parameter $\sol$  and the eigenvalues $\lambda$. In the literature  the a-priori information on the slope parameter $\sol$ such as smoothness is often characterized by considering ellipsoids (see definition below) in $L^2[0,1]$ with respect to a weighted norm $\norm_\bw$ for a pre-specified weight sequence $\bw$.  Moreover, it is usually  assumed that the  sequence $\ev$  of eigenvalues of $\op$ has a polynomial decay (c.f. \cite{HallHorowitz2007} or \cite{CrambesKneipSarda2007}). However, it is well-known that this restriction may exclude several interesting cases, such as an exponential decay. Therefore, we do not  impose a specific  form of a decay. 

It is shown in \cite{Johannes2009} that the estimator $\hsol_m$ given in \eqref{intro:def:est:reg} is optimal in a minimax sense if the parameter $m=m(n)$ is appropriately chosen. Roughly speaking, the introduction of a dimension reduction implies a bias in addition to the classical  variance term which leads the statistician to perform a compromise. The optimal choice of the dimension parameter $m$ requires an a-priori knowledge about the sequences $\bw$ and $\ev$, which is unknown in practice. However, useful elements of this previous work are recalled in Section \ref{sec:bm}.

Our aim in this paper, is to provide a data driven method to select the dimension parameter $m$, in 
such a way that the bias and variance compromise is automatically reached by the resulting estimator. The methodology is inspired by the works of \cite{BBM99}, now extensively described in 
\cite{Massart07} whose results, like ours, are in a non asymptotic setting. By re-writing the estimator $\hsol_m$ as a minimum contrast estimator over the function space $S_m$ $-$ called model $-$ linearly spanned by $\varphi_1, \dots, \varphi_m$, we can propose a model selection device by defining a penalty function. We obtain a selected $\hat m$ in an admissible set of values of $m$. We first define and study in Section \ref{sec:theo1:gen}, the resulting estimator $\hat \beta_{\hat m}$ with deterministic penalty and deterministic set of admissible $m$'s: this requires to assume that the degree of ill-posedness of the problem is known. In other words, information are first supposed to be available about the order of the decay of the eigenvalues $\lambda_j$. This study gives the tools to the next and final step: we define in Section 
\ref{maintheosec} a completely data driven estimator, built by using a random penalty function and a random set of admissible dimensions $m$. We can provide a general risk bound for this estimator and show that it can automatically reach the optimal rate of convergence, without requiring any a-priori knowledge. All proofs are gathered in the Appendix section. 

\section{Background to the methodology.}\label{sec:bm}
\subsection{Notations and basic assumptions}
\paragraph{Circular functional linear model.}In this paper we suppose that  the regressor $X$ is $1$-periodic, that is $X(0)=X(1)$, and second order stationary, i.e., there exists a  positive definite covariance function  $c:[-1,1]\to\R$  such that $\cov(X(t),X(s))=c(t-s)$, $s,t\in[0,1]$. Then it is straightforward to see that the covariance function $c(\cdot)$ is $1$-periodic too. In this situation applying the covariance operator $\op$ equals a convolution with the covariance function. Since $c(\cdot)$ is $1$-periodic it is easily seen that due to the classical convolution theorem, the eigenfunctions of the covariance operator $\op$ are given by the trigonometric basis
\begin{equation*}
\ef_{1}(s) :\equiv1, \;\ef_{2k}(s):=\sqrt{2}\cos(2\pi k s),\; \ef_{2k+1}(s):=\sqrt{2}\sin(2\pi k s),s\in[0,1],\; k\geq1\end{equation*}
and the corresponding eigenvalues satisfy
\begin{equation*}
\ev_{1}=\int_{0}^1c(s)ds,\quad \ev_{2k}=\ev_{2k+1}=\int_{0}^1\cos(2\pi k s) c(s)ds,\; k\geq1. \end{equation*}
Notice that the eigenfunctions are known to the statistician and only the eigenvalues depend on the unknown covariance function $c(\cdot)$, i.e., have to be estimated.

\paragraph{Moment assumptions.}The results  derived below involve  additional  conditions on the moments of the  random function $X$ and the error term $\epsilon$, which we formalize now. Let $\cX$ be the set of all centered $1$-periodic and second order stationary random functions $X\in L^2[0,1]$ with finite second moment, i.e., $\Ex\normV{X}^2<\infty$, and strictly positive covariance operator $\op$. If $\ev:=(\ev_j)_{j\geq1}$ denotes the sequence of eigenvalues associated to $\op$, then given  $X\in\cX$  the random variables $\{[X]_j/\sqrt{\ev_j}, j\in\N\}$  are  centered with variance one.  Here and subsequently, we denote by $\cX^{k}_{\eta}$, $k\in\N$, $\eta\geq1$, the subset of $\cX$  containing only random functions $X$ such that  the $k$-th moment of the corresponding random variables   $[X]_j/\sqrt{\ev_j},$ $j\in\N$ 
are  uniformly bounded, that is
\begin{equation*}
\cX^{k}_{\eta}:=\Bigl\{ X\in\cX\;\text{ such that }\quad\sup_{j\in\N} \Ex\Bigl|[X]_j/\sqrt{\ev_j}\Bigr|^k \leq \eta \Bigr\}.
\end{equation*}
It is worth noting that in case $X\in\cX$ is a Gaussian random function the corresponding random variables   $[X]_j/\sqrt{\ev_j},$ $j\in\N$, are Gaussian with mean zero and variance one. Hence, if $\eta\geq 3$ then  any Gaussian random function $X\in \cX$ belongs also to $\cX^{k}_{\eta}$ for each $k\in\N$.

\paragraph{Minimal regularity conditions.} Given a  strictly positive sequence of  weights $w:=(w_j)_{j\geqslant1}$,  denote by $\cF_{w}^c$ the ellipsoid  with radius  $c>0$, that is, \begin{equation*}
 \cF_{w}^c := \Bigl\{f\in L^2[0,1]: \sum_{j=1}^\infty w_j |\skalarV{f,\ef_j}|^2=:\normV{f}_{w}^2\leq c\Bigr\}.
\end{equation*}
Furthermore, let $\cF_{w}:=\{f\in L^2[0,1]: \normV{f}_{w}^2<\infty\}$ and $\skalarV{f,g}_{w}:= \sum_{j=1}^\infty w_j \skalarV{f,\ef_j}\skalarV{\ef_j,g}$. Note that this weighted inner product induces the weighted norm $\norm_{w}$.

Here and subsequently, given strictly positive sequences of weights $\bw:=(\bw_j)_{j\geqslant1}$ and $\hw:=(\hw_j)_{j\geqslant1}$  we shall  measure the performance of any estimator $\widehat{\beta}$ by its maximal $\cF_{\hw}$-risk over the ellipsoid $\cF_\bw^\br$ with radius $\br>0$, that is $\sup_{\sol\in\cF_\bw^\br}\Ex\|\widehat{\sol}-\sol\|_{\hw}^2$. We do not specify the sequences of weights $\bw$ and $\hw$, but  impose from now on the following minimal regularity conditions.
\begin{assumption}\label{ass:reg} Let $\hw:=(\hw_j)_{j\geqslant 1}$ and $\bw:=(\bw_j)_{j\geqslant 1}$  be positive sequences of weights  with $\hw_1= 1$ and $\bw_1= 1$  such that $(1/\bw_j)_{j\geq 1}$ and   $(\hw_j/\bw_j)_{j\geqslant 1}$ are non increasing zero-sequences.\end{assumption}
Note that under Assumption \ref{ass:reg} the ellipsoid  $\cF_\bw^\br$ is a subset of $\cF_\hw^\br$, and hence the $\cF_{\hw}$-risk a well-defined risk for $\beta$. Roughly speaking, if $\cF_\bw^\br$ describes $p$-times differentiable functions, then the Assumption \ref{ass:reg} ensures that the $\cF_{\hw}$-risk involves maximal $s< p$ derivatives.

\subsection{Minimax optimal estimation.} 
The objective of the paper is to construct an estimator which attains the minimal rate of convergence of the maximal $\cF_{\hw}$-risk over the ellipsoid $\cF_\bw^\br$  for wide range of sequences $\bw$ and $\hw$ satisfying Assumption \ref{ass:reg}, without using an a-priori knowledge of neither $\bw$ nor $\br$.  Therefore, let us first recall a lower bound which can be found in \cite{Johannes2009}.  Let $\mstar:=(\mstarn)\in\N$ for some $\triangle\geq 1$  be chosen such that
\begin{equation*}
1/\triangle\leq \frac{\bw_{\mstarn}}{n\,\hw_{\mstarn}}\sum_{j=1}^{\mstarn} \frac{\hw_{j}}{\lambda_j} \leq \triangle,
\end{equation*} 
i.e. $(1/n)\sum_{j=1}^{\mstarn} \hw_{j}/\lambda_j$ and $\omega_{m_n^*}/\gamma_{m_n^*}$ have the same orders. 
 
Given an i.i.d. $n$-sample of $(Y,X)$ obeying \eqref{intro:e1} with $\sigma>0$ and $X\in\cX$ with associated sequence of eigenvalues  $\ev$, we have then for any estimator $\breve\sol$  that
\begin{equation}\label{res:lower}  \sup_{\sol\in \cF_\bw^\br} \left\{ \Ex\normV{\breve{\sol}-\sol}^2_\hw\right\}\geqslant  \frac{1}{4\,\triangle}\min \left( \frac{\sigma^2}{2},  \frac{\br}{\triangle}\right) \,\max(\hw_{\mstarn}/\bw_{\mstarn},1/n)\quad\mbox{ for all }n\geq1.
\end{equation}
On the other hand  consider the estimator $\hsol_m$ defined in \eqref{intro:def:est:reg} with dimension parameter $m=\mstarn$. If in  addition  $X\in\cX^{16}_{\xi}$, then  it is shown in \cite{Johannes2009} that there exists a numerical constant $C>0$ such that
\begin{gather*}
\sup_{\sol \in \cF_\bw^\br}\left\{\Ex\normV{\hsol_{m_n^*}- \sol}^2_\hw\right\}\leqslant C \, \triangle^3\,\xi\,[ \br\,\Ex\normV{X}^2+\sigma^2]\, \max( \hw_{\mstarn}/\bw_{\mstarn},1/n).
\end{gather*}
Therefore, the minimax-optimal rate of convergence is of order $O(\max( \hw_{\mstarn}/\bw_{\mstarn},1/n))$. As a consequence, the orthogonal series estimator $\hsol_{\mstarn}$ attains this optimal rate  and hence is minimax-optimal.  However, the definition of the dimension parameter $\mstarn$ used to construct the estimator involves an a-priori knowledge of the  sequences $\bw$, $\hw$ and $\ev$. Throughout the paper our aim is to construct a data-driven choice of the dimension parameter  not requiring this a-priori knowledge and automatically  attaining the optimal rate of convergence.

\subsection{Example of rates}\label{illustre}
We compute in this section the rates that we can obtain in three
configurations for the sequences $\bw, \hw $ and $ \lambda$. These
cases will be referred to in the following. In all three cases, we take the sequence $\hw$ with $\hw_j=j^{2s}$, $j\geq1$, for $s\in {\mathbb R}$. 
\paragraph{Case [P-P] Polynomial-Polynomial.} Consider sequences $\bw$ and $\ev$ with $\bw_j=j^{2p}$, $j\geq1$, for $p>\max(0,s)$,
and $\ev_j\asymp j^{-2a}$, $j\geq 1$, for $a>1/2$ respectively, where the notation $u_j\asymp v_j$, $j\geq1$, means that there exists a constant $d>0$ such that $u_j/d \leq  v_j \leq du_j$ for all $j\geq 1$. Then it is easily seen that $(\mstarn)^{2(s-p)}=\frac{\hw_{\mstarn}}{\bw_{\mstarn}}
\asymp \sum_{j=1}^{\mstarn} \frac{\hw_{j}}{n\lw_j} \asymp n^{-1} \sum_{j=1}^{\mstarn} j^{2s+2a}$ and hence $m_n^*\asymp n^{1/(2p+2a+1)}$ if $2s+2a+1>0$, $m_n^*\asymp n^{1/[2(p-s)]}$ if
$2s+2a+1<0$ and $m_n^*\asymp(n/\log(n))^{1/[2(p-s)]}$ if $2a+2s+1=0$. Finally, 
the optimal rate attained by the estimator is 
$\max(n^{-(2p-2s)/(2a+2p+1)},n^{-1})$, if $2s+2a+1\neq 0$ (and
$\log(n)/n$ if $2s+2a+1=0$).  Observe that an increasing value of $a$ leads to a slower optimal rate of convergence. Therefore, the parameter $a$ is called degree of ill-posedness (c.f. \cite{Natterer84}).
\begin{rem} Obviously the rate is parametric if $2a+2s+1<0$. The case $0\leq s <p$ can be
interpreted as the $L^2$-risk of an estimator of the $s$-th derivative of the slope parameter $\sol$. On the other hand
the case, $s=-a$, corresponds to the mean-prediction error (c.f. \cite{CardotJohannes2007}).\hfill$\square$
\end{rem}
\paragraph{Case [E-P] Exponential-Polynomial.} Consider sequences $\bw$ and $\ev$ with 
$\bw_j=\exp (j^{2p})$, $j\geq1$, for $p>0$, and  (as previously) $\ev_j\asymp j^{-2a}$, $j\geq 1$, for $a>1/2$ respectively.
Then $\mstarn$ is such that  $\exp(-(\mstarn)^{2p})(\mstarn)^{2s}=\frac{\hw_{\mstarn}}{\bw_{\mstarn}}
\asymp \sum_{j=1}^{\mstarn} \frac{\hw_{j}}{n\lw_j} \asymp n^{-1} \sum_{j=1}^{\mstarn} j^{2s+2a}$. In case $2a+2s+1>0$ this is equivalent to $\exp(-(\mstarn)^{2p})\asymp (\mstarn)^{2a+1}n^{-1}$ and hence $\mstarn\asymp (\log n -\frac{2a+1}{2p}\log(\log n))^{1/(2p)}.$ Thereby,  $n^{-1} (\log n)^{(2a+1+2s)/(2p)}$ is the optimal rate attained by the estimator.
Furthermore, if $2a+2s+1<0$, then $m_n^*\asymp(\log(n)+(s/p)\log(\log(n)))^{1/(2p)}$ and the rate is parametric, while if $2a+2s+1=0$, the rate is of order $\log(\log(n))/n$.

\paragraph{Case [P-E] Polynomial-Exponential.} Consider sequences $\bw$ and $\ev$ with $\bw_j=j^{2p}$, $j\geq1$, for $p>\max(0,s)$,  and  $\ev_j\asymp \exp(-j^{2a})$, $j\geq 1$, for $a>0$ respectively. Then $(\mstarn)^{2(s-p)}=\frac{\hw_{\mstarn}}{\bw_{\mstarn}}
\asymp \sum_{j=1}^{\mstarn} \frac{\hw_{j}}{n\lw_j} \asymp n^{-1} \sum_{j=1}^m j^{2s}\exp(j^{2a})$ and hence $\mstarn\asymp (\log n -\frac{2p+(2a-1)_{\vee0}}{2a}\log(\log n))^{1/(2a)}$ with $(q)_{\vee0}:=\max(q,0)$. Thereby, $(\log n)^{-(p-s)/a}$ is the optimal rate attained by the estimator. The parameter $a$ reflects again the  degree of ill-posedness since an increasing value of $a$ leads also here to a slower optimal rate of convergence.

\section{A model selection approach: known degree of ill-posedness}\label{sec:theo1:gen}
In the previous section, we have recalled an estimation procedure that attains the optimal rate of convergence in case the slope parameter belongs to some ellipsoid $\cF_\bw^\br$ and its accuracy is measured by a $\cF_\hw$-risk.  In this section, we suppose that there exists an a-priori knowledge concerning the degree of ill-posedness, that is the asymptotic behavior of the sequence of eigenvalues $\ev$ is known. The objective is the construction of an adaptive estimator which depends neither on the sequence of weights $\bw$ nor on the radius $\br$ but still attains the optimal rate over the ellipsoid $\cF_\bw^\br$.
In this section, we use the following assumption.
\begin{assumption}\label{ass:reg:ext1}Let $\ev:=(\ev_j)_{j\geq1}$ denote the sequence of eigenvalues associated to
the regressor $X$ and let $\hw:=(\hw_j)_{j\geq1}$  be a sequence satisfying Assumption \ref{ass:reg} such that 
\begin{itemize}
 \item[(i)] there exist non decreasing sequences $\delta:=\delta(\ev,\hw):=(\delta_m(\ev,\hw))_{m\geq1}$ and $\Delta:=\Delta(\ev,\hw):=(\Delta_m(\ev,\hw))_{m\geq1}$ with  $\delta_m\geq \sum_{j=1}^m\hw_j/\ev_j$ and $\Delta_m\geq \max_{1\leq j\leq m}\hw_j/\ev_j$ for all $m\geq 1$ such  that for some  $\Sigma>0$,%
\begin{equation}\label{sumcondi} \sum_{m\geq 1}\Delta_m\exp(-\frac{\delta_m}{6\Delta_m})\leq \Sigma.\end{equation}
\item[(ii)] the sequence  $M:=(M_n)_{n\geq1}$  given by $ M_n:=\argmax_{1\leq M\leq n}\{\delta_{M}\leq \delta_1 n {(\hw_M)_{\wedge1}}\}$, $n\geq 1$, with $(q)_{\wedge1}:=\min(q,1)$, satisfies \begin{equation}\label{ass:reg:ext1:1}  \min_{1\leq j\leq  M_n} \ev_{j} \geq 2/n\qquad\mbox{ for all $n\geq1$.}\end{equation}
\end{itemize}
\end{assumption}
It is worth to note that both  sequences $\delta$ and $M$ depend on the eigenvalues $\ev$.
\subsection{Definition of the estimator.} Consider the orthogonal series estimator $\hsol_{m}$ defined in
\eqref{intro:def:est:reg}. In what follows we construct an adaptive procedure to choose the
dimension parameter $m$ based on a model selection approach.
Therefore, let $\hPhi_{u}= \sum_{j\geq 1}\hev_j^{-1}\1\{\hev_j\geq
1/n\}[u]_j\ef_j$ for $u\in L^2[0,1]$ with Fourier coefficients
$[u]_j:=\skalarV{u,\ef_j}$. Then we consider the contrast
\begin{equation}\label{sec:gen:con:1}
\ct(t) :=   \normV{t}^2_\hw - 2 {\skalarV{ t, \hPhi_{\hg}}_\hw}.
\end{equation}
Define $\cS_m:={\rm span}\{\ef_1,\dotsc,\ef_m\}$. Obviously for all $t\in\cS_m$ it follows that  $\skalarV{t,\hPhi_{\hg}}_\hw=\skalarV{t,\hsol_m}_\hw$ and hence $\ct(t) = \normV{t-\hsol_m}^2_\hw - \normV{\hsol_m}^2_\hw$. Therefore, we have  for all $m\geq1$
\begin{equation*}
\arg\min_{t\in\cS_m}  \ct(t) = \hsol_{m}.
\end{equation*}
Let $X\in\cX^{4}_{\eta}$  and  $\Ex |Y/\sigma_Y|^{4} \leq \eta$ with $\sigma_Y^2:=\Var(Y)$. Under Assumption  \ref{ass:reg:ext1}, we consider the penalty function
\begin{equation*}
\pen(m) := 192\sigma_Y^2 \eta \frac{\delta_m}{n}.
\end{equation*}
The adaptive estimator $\hsol_{\whm}$ is obtained from \eqref{intro:def:est:reg} by choosing the dimension parameter
\begin{equation}\label{sec:gen:ada:1}
\whm := \argmin_{1\leq m\leq M_n}\left\{\ct(\hsol_{m}) + \pen(m)   \right\}.
\end{equation}
Note that we can compute $$\ct(\hsol_{m})=-\sum_{j=1}^m \omega_j\frac{[\widehat g]_j^2}{\widehat
\lambda_j^2}\1\{\hev_{j}\geqslant 1/n\}.$$
\begin{rem}\label{rem:gen:ada:1}Throughout the paper we ignore that also the value $\sigma_Y^2$
and $\eta$ are unknown in practice. Obviously $\sigma_Y^2$ can be
estimated straightforwardly by its empirical counterpart.  An
estimator of the value $\eta$ is not a trivial task. However, if  in
addition the regressor $X$ and the error term $\epsilon$ are
Gaussian, then $Y\sim \cN(0,\sigma_Y^2)$ and hence $\eta=3$ is
a-priori known. We may take an other point of view if we chose a-priori a sufficiently large $\eta\geq 3$ (the Gaussian case is included) then the following assertions apply as long as the unknown data
generating process satisfies the conditions $X\in\cX^{4}_{\eta}$  and  $\Ex |Y/\sigma_Y|^{4} \leq \eta$.\hfill$\square$\end{rem}

\subsection{An upper bound.} We derive first an upper bound of the
adaptive estimator $\hsol_{\whm}$ by assuming an a-priori knowledge
of appropriate sequences $\delta$  and $M$ which are used in the construction of the penalty and the admissible set of values of $m$.
\begin{theo}\label{sec:gen:th1}Assume an $n$-sample of $(Y,X)$ satisfying \eqref{intro:e1}. Let $\Ex |Y/\sigma_Y|^{4} \leq \eta$ and $X\in\cX^{4}_{\eta}$ be  $1$-periodic and second order stationary with associated  eigenvalues $\ev$.

Suppose that the sequences $\bw$ and $\hw$    satisfy Assumption \ref{ass:reg}. Let $\delta$,  $\triangle$  and $M$ be sequences
satisfying Assumption \ref{ass:reg:ext1} for some constant $\Sigma$. Consider
the estimator  $\hsol_{\whm}$   defined in \eqref{intro:def:est:reg}
with  $\whm$ given by \eqref{sec:gen:ada:1}.  If in addition $X\in\cX^{24}_{\xi}$ and  $\Ex
|Y/\sigma_Y|^{24} \leq \xi$, then there exists a
numerical constant $C$ such that for all $n\geq1$ and $1\leq m\leq
M_n$,  we have
\begin{multline*}
\sup_{\sol \in \cF_\bw^\br}\left\{\Ex\normV{\hsol_{\whm}-
\sol}^2_\hw\right\}\leqslant C \Bigl\{ \frac{\hw_{m}}{\bw_{m}}
\,\br+  \frac{\delta_{m}}{n} \,(\br \Ex\normV{X}^2 +\sigma^2)\eta \,\Bigr\}\\\hfill  +
\frac{K}{n}\,  (\br \Ex\normV{X}^2 +\sigma^2)\, [ \delta_1 + \br][ 1 +
(\Ex\normV{X}^2)^2],
\end{multline*}
where $K=K(\Sigma, \eta,\xi,\delta_1)$ is a constant depending on $\Sigma,  \eta,\xi$ and $\delta_1$ only.
\end{theo}

It is worth noting, that in the last assertion we do not impose a
complete knowledge of the sequence of eigenvalues $\ev$ associated
to the regressor $X$. In the next Corollary we state the upper bound when balancing the terms depending on $m$, which is obviously a trivial consequence of Theorem \ref{sec:gen:th1}.
\begin{coro}\label{sec:gen:prop1:cor}
Let the assumptions of Theorem \ref{sec:gen:th1} be satisfied. If in addition the sequence $\mdag:=(\mdagn)_{n\geq1}$  is chosen such that
${\bw_{\mdagn}\delta_{\mdagn}}/({n\,\hw_{\mdagn}})\asymp 1$, $n\geq 1$, then we have
\begin{gather*}
\sup_{\sol \in \cF_\bw^\br}\left\{\Ex\normV{\hsol_{\whm}- \sol}^2_\hw\right\}=O\Bigl( \max({\hw_{\mdagn}}/{\bw_{\mdagn}}, {1}/{n})\Bigr) \mbox{ as }n\to\infty.
\end{gather*}
\end{coro}
\begin{rem} Comparing the last assertion with the lower bound given in \eqref{res:lower}, we see that the adaptive estimator attains the optimal rate of convergence, as long as\linebreak $\sup_{n\geq1}{\hw_{\mdagn}\bw_{\mstarn}}/({\bw_{\mdagn}\hw_{\mstarn}})<\infty$. Obviously a sufficient condition is given if the sequence $\delta$ satisfies  in addition   $\sup_{m\geq1}\delta_m/(\sum_{j=1}^m\hw_j/\ev_j)<\infty$. The polynomial case below provides an example. However, this condition is not necessary as can be seen in the exponential case. \hfill$\square$\end{rem}

\subsection{Convergence rate of the theoretical adaptive estimator.}
We described in Section \ref{illustre} three different cases where
we could choose the model $m$ such  that the resulting estimator
reaches the optimal minimax rate. The following result shows that, in case of known degree of ill-posedness, we can
propose choices of sequences $\delta$, $\Delta$ and $M$ such that the penalized
estimator automatically attains the optimal rate.

\begin{prop}\label{ratesknownlambda1} In cases {\bf [P-P]} and {\bf [E-P]} with $2a+2s+1>0$, let
$\delta_m\asymp m^{2a+2s+1}$,  $\Delta_m\asymp
m^{(2a+2s)_{\vee0}}$ and  $M_n\asymp n^{1/(2a+1+(2s)_{\vee 0})}$with $(q)_{\vee0}:=\max(q,0)$.
While in case {\bf [P-E]}, choose $\delta_m\asymp m^{2a+1+(2s)_{\vee0}}\exp(m^{2a})$,   $\Delta_m\asymp m^{(2s)_{\vee0}}\exp(m^{2a})$ and $M_n\asymp(\log n/(\log n)^{(2a+1+(2s)_{\vee0})/(2a)})^{1/(2a)}$.

Then Assumption \ref{ass:reg:ext1} is fulfilled and, under the
additional assumptions of Theorem \ref{sec:gen:th1}, the adaptive
estimator $\widehat\beta_{\widehat m}$ reaches the optimal rate.
\end{prop}

In cases {\bf [P-P]} and {\bf [E-P]}, if $2a+2s+1< 0$,  then the sequence $\delta$
can be taken of order 1. The collection of models must be reduced to
$\{[\sqrt{n}], \dots, n\}$ since $M_n$ can be taken equal to $n$. It
appears then that the rate is parametric in this case. In fact, no
model selection is necessary in this case, a large $m$ ($m=n$ for
instance) can be chosen.\\

Now, we have in mind to prepare the case where  the degree of ill-posedness of the $\ev_j$'s,  and  more precisely  
$\delta_m$ and $M_n$, are unknown. We propose
hereafter a more intrinsic choice of $\delta_m$, which does not
require anything but the $\lambda_j$'s (which can be estimated). In this spirit, we can prove the following assertion.

\begin{prop}\label{ratesknownlambda2}
In cases {\bf [P-P]} and {\bf [E-P]} with $a+s\geq 0$ or in case
{\bf [P-E]}, choose $\Delta_m:=\max_{1\leq j\leq m}\hw_j/\ev_j$,
$\kappa_m:=\max_{1\leq j\leq m}(\hw_j)_{\vee1}/\ev_j$ with $(q)_{\vee1}:=\max(q,1)$ and
\begin{equation}\label{deltamref} \delta_m :=  m
\Delta_m\Bigl|\frac{\log (\kappa_m\vee (m+2))}{\log(m+2)}\Bigr|.
\end{equation}
Then Assumption \ref{ass:reg:ext1} is fulfilled and, under the
additional assumptions of Theorem \ref{sec:gen:th1}, the adaptive
estimator $\widehat\beta_{\widehat m}$ reaches the optimal rate.
\end{prop}

\section{A model selection approach: unknown degree of ill-posedness}\label{maintheosec}
In this section, the objective is the
construction of a fully adaptive estimator which does not depend on the
sequence $\bw$ and $\lambda$. Nevertheless the resulting estimator still attains the optimal rate in
case the slope parameter belongs to some ellipsoid $\cF_\bw^\br$ and
the sequence of eigenvalues $\lambda$ associated to the covariance operator of
$X$ has a given (unknown) rate of decrease.

The configuration given in Proposition \ref{ratesknownlambda2} is
now the right reference and the choice that the estimator is going
to mimic.  In particular, it is easily seen that there exists always a constant $\Sigma>0$  such that the sequences $\delta$ and $\triangle$  given in Proposition \ref{ratesknownlambda2} satisfy  Assumption \ref{ass:reg:ext1} (i).
Observe that in this situation we have
\begin{eqnarray*}\Delta_m\exp(-\frac{\delta_m}{6\Delta_m}) &= &\Delta_m\exp(-\frac{m }{6} \frac{\log (\kappa_m\vee (m+2))}{\log(m+2)})\\ &\leq & (\kappa_m\vee (m+2))\exp(-\frac{m }{6} \frac{\log (\kappa_m\vee (m+2))}{\log(m+2)})\\ &\leq & \exp\Bigl(-m\Bigl[\frac{1}{6}-\frac{\log(m+2)}{m}\Bigr] \frac{\log (\kappa_m\vee (m+2))}{\log(m+2)}\Bigr) 
\end{eqnarray*}
where the last term is obviously  summable.

\begin{assumption}\label{ass:reg:ext2} Let $\ev$ denote the sequence of eigenvalues associated to
the regressor $X$, let $\delta$ and $\triangle$ be  the sequences defined in Proposition \ref{ratesknownlambda2} and let $\bw$ and $\hw$  be sequences satisfying Assumption \ref{ass:reg} such that
\begin{itemize}
\item[(i)]the sequence  $M:=(M_n)_{n\geq1}$ given in Assumption \ref{ass:reg:ext1} satisfies in addition to \eqref{ass:reg:ext1:1} also
\begin{equation*}
 \frac{\log n}{2n} \geq \max_{m> M_n} \frac{\ev_{m}}{m(\hw_m)_{\vee1}} \qquad\mbox{ for all $n\geq1$;}\end{equation*}
\item[(ii)]the sequence $\mdag:=(\mdagn)_{n\geq1}$ given by $1/\uc\leq {\bw_{\mdagn}\delta_{\mdagn}}/({n\,\hw_{\mdagn}})\leq \uc$ for all $n\geq 1$ and some $\uc\geq 1$ satisfies
\begin{equation*}\min_{1\leq m\leq\mdagn} \frac{\ev_{m}}{m(\hw_m)_{\vee1}}\geq 2  (\log n)/n\qquad\mbox{ for all $n\geq1$;} \end{equation*}
\item[(iii)]the  sequence $N:=(N_n)_{n\geq1}$ given by $N_n:=\argmax_{1\leq N\leq n}\{ \max\limits_{1\leq j\leq N} \hw_{j}/n\leq 1\}$, $n\geq 1$, satisfies 
\begin{equation*}  M_n\leq N_n\leq n  \qquad\mbox{ for all $n\geq1$.} \end{equation*}
\end{itemize}
\end{assumption}
\begin{rem}The last assumption is  technical
but satisfied in the interesting case. Note that $(i)$ and $(ii)$
together imply $\mdagn\leq M_n$ for all $n\geq1$. The condition $(iii)$
is rather weak, observe that the sequence $\hw$ is a-priori known
and thus also the sequence of upper bounds $N$. In particular, recall that in case $\hw\equiv 1$ the $\cF_\hw$-risk corresponds to the $L^2$-risk. If $\hw_m\leq 1$ for all $m\geq1$, then  $\cF_\hw$-risk is weaker than the $L^2$-risk and  $N_n=n$. Only if the $\cF_\hw$-risk is stronger than the
$L^2$-risk, that is $\hw$ is monotonically increasing, we choose
$N_n$ such that $\hw_{N_n} \asymp n$. Then it is not hard to see
that in these situations $(iii)$ is satisfied at least for sufficiently large $n$. \hfill$\square$\end{rem}
\subsection{Definition of the estimator}We follow the model
selection approach presented in the last section. Define 
\begin{equation*}
\widehat \Delta_m:=\max_{1\leq j\leq  m}\frac{\hw_j}{\hev_j}
\1_{\{\hev_j\geq 1/n\}}\quad\mbox{ and }\quad\widehat\kappa_m:=\max_{1\leq j\leq
m}\frac{(\hw_j)_{\vee 1}}{\hev_j} \1_{\{\hev_j\geq 1/n\}}.
\end{equation*}
We shall refer to $\delta_m$ as
defined in (\ref{deltamref}) and consider its estimator given by\begin{equation*}\widehat \delta_m :=  m  \widehat\Delta_m\Bigl|\frac{\log (\widehat\kappa_m\vee (m+2))}{\log(
m+2)}\Bigr|.\end{equation*} 
If $X\in\cX^{4}_{\eta}$ and $\Ex |Y/\sigma_Y|^{4} \leq \eta$, then we define a random penalty function
\begin{equation*}
\hpen(m)= 1920 \sigma_Y^2 \eta\; \frac{\widehat \delta_m }{n}.
\end{equation*}
Moreover, we consider a random upper bound for the collection of models given by
\begin{equation}\label{sec:gen:rub}\hM_n:= \argmax_{1\leq M\leq N_n} \Bigl\{ \frac{\hev_{M}}{M(\hw_M)_{\vee1}} \geq {(\log n)}/{n}\Bigr\}.\end{equation}
The adaptive estimator $\hsol_{\whm}$ is obtained from \eqref{intro:def:est:reg} by choosing the dimension parameter
\begin{equation}\label{sec:gen:ada:2}
\whm := \argmin_{1\leq m\leq \hM_n}\left\{\ct(\hsol_{m}) + \hpen(m)   \right\}
\end{equation}
We shall emphasize that the proposed estimator does not depend on an a-priori knowledge of
neither the sequence $\bw$ nor the sequence $\ev$. 

\subsection{An upper bound.} In the next assertion we provide an upper bound of the fully 
adaptive estimator $\hsol_{\whm}$ by assuming that the sequences $\ev$,
$\hw$ and $\bw$ satisfy Assumption \ref{ass:reg:ext2}.
\begin{theo}\label{sec:gen:th2}Assume an $n$-sample of $(Y,X)$ satisfying \eqref{intro:e1}.  Suppose that $\Ex |Y/\sigma_Y|^{4} \leq \eta$ and that $X\in\cX^{4}_{\eta}$ is $1$-periodic and second order stationary. Let Assumption \ref{ass:reg:ext2} be satisfied. Consider the estimator $\hsol_{\whm}$ defined in \eqref{intro:def:est:reg} with ${\whm}$ given by \eqref{sec:gen:ada:2}. If in  addition $X\in\cX^{28}_{\xi}$ and $\Ex |Y/\sigma_Y|^{28} \leq \xi$,  then there exists a numerical constant $C>0$  such that for all $n\geq1$%
\begin{multline*}
\sup_{\sol \in \cF_\bw^\br}\left\{\Ex\normV{\hsol_{\whm}- \sol}^2_\hw\right\}\leqslant  C\,\frac{\hw_{\mdagn}}{\bw_{\mdagn}}(\br + \uc\eta[\br \Ex\normV{X}^2 +\sigma^2])\\ +     \frac{K}{n}\,  [\br \Ex\normV{X}^2 +\sigma^2]\, [ 1+\delta_1 + \br][ 1 + 
(\Ex\normV{X}^2)^2],
\end{multline*}
where $\mdagn$ and $\uc$ are  defined in Assumption \ref{ass:reg:ext2},  $K=K(\Sigma, \eta,\xi,\delta_1)$  is a constant only depending on $ \eta,\xi,\delta_1$ and $\Sigma$  such that the sequences $\delta$ and $\triangle$  given in Proposition \ref{ratesknownlambda2} satisfy Assumption \ref{ass:reg:ext1}.
\end{theo}
\begin{rem}Comparing the last assertion with Theorem  \ref{sec:gen:th1}, we see
that under Assumption \ref{ass:reg:ext2} the proposed adaptive
estimator obtains the same rate as in case of known degree of
ill-posedness. We only have to impose in addition slightly stronger
moment conditions.\hfill$\square$\end{rem}
It is easily verified that in all the examples discussed above the fully adaptive estimator attains the optimal rate, which is summarized in the next assertion. 
\begin{coro}\label{optimadapt}
In cases {\bf [P-P]} and {\bf [E-P]} with $a+s\geq 0$ or in case
{\bf [P-E]}, Assumption \ref{ass:reg:ext2} is fulfilled and, under
the additional assumptions of Theorem \ref{sec:gen:th2}, the fully 
adaptive estimator $\widehat\beta_{{\whm}}$ with ${\whm}$ given by \eqref{sec:gen:ada:2} reaches
the optimal rate.
\end{coro}

\paragraph{Conclusion.}

Assuming a circular functional linear model we derive in this paper a fully adaptive estimator of the slope function $\beta$ or its derivatives,  which attains the  minimax optimal rate of convergence. It is worth to note, that in this paper not only the penalty is chosen randomly but also the collection of models. In this way the proposed estimator is adaptive also with respect to the degree of ill-posedness of the underlying inverse problem.  We can thereby face both, the  mildly and the severely ill-posed case.

It is not clear that the ideas in this paper can be straightforwardly adapted to treat the case of noncircular functional models. We are currently exploring this issue.

\appendix
\section{Appendix}
\subsection{Proof of Theorem \ref{sec:gen:th1}}
We begin by defining and recalling notations to be used in the proof. Given $u\in L^2[0,1]$ we denote by $[u]$ the infinite vector of Fourier coefficients $[u]_j:=\skalarV{u,\ef_j}$. In particular
we use the notations 
\begin{multline*}
[X_{i}]_j= \skalarV{X_i,\ef_j}, \quad [\sol]_j=\skalarV{\sol,\ef_j},\quad \sigma_Y^2=\Var(Y),\\
\hsol_m= \sum_{j=1}^m\hev_j^{-1}\1\{\hev_j\geq 1/n\}[\hg]_j \ef_j, \quad\tsol_m:= \sum_{j=1}^m
\ev_j^{-1}[\hg]_j \ef_j,\quad{\sol}_m:=\sum_{j=1}^m [\sol]_{j}  \ef_{j}, \\
\hPhi_{u}= \sum_{j\geq 1}\hev_j^{-1}\1\{\hev_j\geq 1/n\}[u]_j\ef_j,\quad \tPhi_{u}:= \sum_{j\geq 1}\lambda_j^{-1}[u]_j\phi_j.
\hfill\end{multline*} 
Given $m\geq1$  we have then for all $t\in\cS_m={\rm span}\{\ef_1,\dotsc,\ef_m\}$ 
 \begin{multline}\label{app:proofs:gen:1:def:2}
\skalarV{ t, \sol}_\hw=  \sum_{j=1}^m \hw_j[t]_j[\sol]_j= \sum_{j=1}^m\frac{\hw_j[t]_j[g]_j}{\ev_j}= \skalarV{ t,
\tPhi_{g}}_\hw,\\
 \skalarV{ t, \tsol_m}_\hw=\skalarV{ t,\tPhi_{\hg}}_\hw=  \frac{1}{n} \sum_{i=1}^n Y_i 
\skalarV{t,\tPhi_{X_i}}_\hw=\frac{1}{n} \sum_{i=1}^n Y_i \sum_{j=1}^m
\frac{\hw_j}{\ev_j}[X_i]_j[t]_j,
\\
\skalarV{ t, \hsol_m}_\hw=  \skalarV{ t,
\hPhi_{\hg}}_\hw=\frac{1}{n} \sum_{i=1}^n Y_i 
\skalarV{t,\hPhi_{X_i}}_\hw=\frac{1}{n} \sum_{i=1}^n Y_i \sum_{j=1}^m\frac{\hw_j}{\hev_j}\1\{\hev_j\geq 1/n\}[X_i]_j[t]_j.\hfill
\end{multline}
Furthermore,  define the event \begin{equation*}
                               \Omega_{Y,X}:=\{ |Y/\sigma_Y| \leq n^{1/6}, |[X]_j/\sqrt{\lambda_j}|\leq n^{1/6}, 1\leq j\leq M_n\}
                               \end{equation*}and denote its complement by $\Omega_{Y,X}^c$. Then consider the functions $\widehat h$ and $\widehat f$ with Fourier coefficients given by  
\begin{multline*}
[\widehat h]_j:= \frac{1}{n}\sum_{i=1}^n \{ Y_i[X_i]_j \1\Omega_{Y_i,X_i} - \Ex (Y_i[X_i]_j \1\Omega_{Y_i,X_i})\},\\
[\widehat f]_j:= \frac{1}{n}\sum_{i=1}^n \{ Y_i[X_i]_j \1\Omega_{Y_i,X_i}^c - \Ex (Y_i[X_i]_j \1\Omega_{Y_i,X_i}^c)\}.\hfill\end{multline*} 
Obviously we have $ [\hg]_j-[g]_j= [\widehat h]_j + [\widehat f]_j $ and hence for all $t\in\cS_m$ 
 \begin{eqnarray}\nonumber
 \skalarV{ t,\hPhi_{\hg}-\sol}_\hw &=& \skalarV{ t,\hPhi_{\hg}-\tPhi_{g}}_\hw=  \skalarV{ t,\tPhi_{\hg}-\tPhi_{g}}_\hw +\skalarV{ t,\hPhi_{\hg}-\tPhi_{\hg}}_\hw\\
\label{app:proofs:gen:1:def:4}  &= & \skalarV{ t,\tPhi_{\widehat h}}_\hw  +\skalarV{ t,\tPhi_{\widehat f}}_\hw  +\skalarV{ t,\hPhi_{\hg}-\tPhi_{\hg}}_\hw  .
\end{eqnarray}
We shall prove in the end of this section three technical Lemmas (\ref{app:th1:l1} - \ref{app:th1:l4}) which are used in the following steps of the proof.\\[1ex]

Consider now the contrast $\ct$
 then  by using  \eqref{sec:gen:con:1} and \eqref{sec:gen:ada:1}  it follows that
\begin{equation*}
\ct(\hsol_{\whm}) + \pen(\whm) \leq \ct(\hsol_{m}) + \pen(m) \leq  \ct(\sol_m) + \pen(m),\qquad\forall1\leq m \leq M_n,
\end{equation*} which in particular implies by using the notations given in \eqref{app:proofs:gen:1:def:2}  that 
\begin{multline*} 
\normV{\hsol_{\whm}}_\hw^2-\normV{\sol_m}_\hw^2\leq 2\{\skalarV{\hsol_{\whm}, \hPhi_{\hg}}_\hw   -\skalarV{\sol_{m}, \hPhi_{\hg}}_\hw\} + \pen(m)-\pen(\whm)\\
=2\skalarV{\hsol_{\whm}-\sol_m, \hPhi_{\hg}}_\hw  + \pen(m)-\pen(\whm) .\end{multline*}
Rewriting the last estimate by using \eqref{app:proofs:gen:1:def:4} we conclude that
\begin{align}\nonumber
\normV{\hsol_{\whm}-\sol}^2_\hw&= \normV{\sol-\sol_{m}}^2_\hw + \normV{\hsol_{\whm}}^2_\hw   - \normV{\sol_{m}}^2_\hw  - 2 \skalarV{\hsol_{\whm}-\sol_m,\sol}_\hw\\\nonumber
&\leq \normV{\sol-\sol_{m}}^2_\hw + \pen(m)-\pen(\whm)+ 2\skalarV{\hsol_{\whm}-\sol_m, \hPhi_{\hg}-\sol}_\hw\\\nonumber
&\leq \normV{\sol-\sol_{m}}^2_\hw + \pen(m)-\pen(\whm) \\\label{pr:th1:e2}
&\quad+2\skalarV{ \hsol_{\whm}-\sol_m,\tPhi_{\widehat h}}_\hw  +2\skalarV{ \hsol_{\whm}-\sol_m,\tPhi_{\widehat f}}_\hw  +2\skalarV{ \hsol_{\whm}-\sol_m,\hPhi_{\hg}-\tPhi_{\hg}}_\hw.
\end{align}
Consider the unit ball $\cB_m:=\{f\in \cS_m:\normV{f}_\hw\leq 1\}$  and  let $\whm\vee m:=\max(\whm,m)$.  Combining  for $\tau>0$ and $ f\in \cS_m$ the elementary inequality
\begin{equation*}
2|\skalarV{f,g}_\hw|
\leq 2\normV{f}_\hw \sup_{t\in \cB_m}|\skalarV{t,g}_\hw|\leq \tau \normV{f}^2_\hw+\frac{1}{\tau} \sup_{t\in \cB_m}|\skalarV{t,g}_\hw|^2
\end{equation*} 
with  \eqref{pr:th1:e2} and $\hsol_{\whm}-\sol_m \in \cS_{\whm\vee m} \subset \cS_{M_n}$ we obtain 
\begin{multline*}
\normV{\hsol_{\whm}-\sol}^2_\omega \leq \normV{\sol-\sol_{m}}^2_\omega + 6 \tau \normV{\hsol_{\whm}-\sol_m}^2_\omega  + \pen(m)-\pen(\whm) \\
+ \frac{2}{\tau}\sup_{ t\in \cB_{\whm\vee m}}  |\skalarV{ t,\tPhi_{\widehat h}}_\hw|^2  +   \frac{2}{\tau}\sup_{ t\in \cB_{M_n} } |\skalarV{ t,\tPhi_{\widehat f}}_\hw|^2 +   {\frac{2}{\tau}\sup_{ t\in \cB_{M_n} } |\skalarV{ t,\hPhi_{\hg}-\tPhi_{\hg}}_\hw|^2}. 
\end{multline*}
Then, noting that $\pen(m\vee m')\leq \pen(m) +\pen(m')$ and $\normV{\hsol_{\whm}-\sol_m}^2_\omega \leq 2\normV{\hsol_{\whm}-\sol}^2_\omega +2\normV{\sol_m-\sol}^2_\omega$, we get, together for $\tau=1/16$ and $\pen(m)= 192\sigma_Y^2\eta  \delta_{m}/n$ that
\begin{multline}\label{pr:th1:e3}
(1/4) \normV{\hsol_{\whm}-\sol}^2\leq (7/4) \normV{\sol-\sol_{m}}^2 + 32\Bigl(\sup_{ t\in \cB_{\whm\vee m}}  |\skalarV{ t,\tPhi_{\widehat h}}_\hw|^2-(1/{32})\pen(\whm\vee m)\Bigr)_+\\
\hfill+   32 \sup_{ t\in \cB_{M_n} } |\skalarV{ t,\tPhi_{\widehat f}}_\hw|^2+ { 32
\sup_{ t\in \cB_{M_n} } |\skalarV{ t,\hPhi_{\hg}-\tPhi_{\hg}}_\hw|^2} +   \pen(\whm\vee m)  + \pen(m)-\pen(\whm)
\\
\hfill\leq (7/4) \normV{\sol-\sol_{m}}^2 + 32\sum_{m'=1}^{M_n}\Bigl(\sup_{ t\in \cB_{ m'}}  |\skalarV{ t,\tPhi_{\widehat h}}_\hw|^2-6\sigma_Y^2\eta \delta_{m'}/n\Bigr)_+\hfill\\
\hfill+   32 \sup_{ t\in \cB_{M_n} } |\skalarV{ t,\tPhi_{\widehat f}}_\hw|^2+ { 32
\sup_{ t\in \cB_{M_n} } |\skalarV{ t,\hPhi_{\hg}-\tPhi_{\hg}}_\hw|^2} +  2 \pen(m).
\end{multline}
Combining the last bound with \eqref{app:th1:l1:e1}  in Lemma \ref{app:th1:l1},  \eqref{app:th1:l2:e1}  and \eqref{app:th1:l2:e2}  in Lemma \ref{app:th1:l2} we conclude that there exist a numerical constant $C$  and a constant $K(\Sigma, \eta)$ depending on $\Sigma$ and $ \eta$ only, such that for all $n\geq1$ and for all $1\leq m \leq M_n$ we have 
$$
\Ex \normV{\hsol_{\whm}-\sol}^2_\hw\leq 7\normV{\sol-\sol_{m}}^2_\hw +8\pen(m)+ \frac{1}n  [C\xi (\sigma^2_Y\delta_1 + \normV{\sol}_\hw^2\}\{ 1 + (\Ex\normV{X}^2)^2\} + \sigma^2_Y K(\Sigma(6),\eta)].
$$ 
Since  $(\hw/\bw)$ is monotonically non increasing we obtain in case $\sol\in\cF_\bw^\br$ that  $\normV{\sol}^2_\hw\leq \br$ and $\normV{\sol-\sol_{m}}^2_\hw\leq (\hw_m/\bw_{m})\br$. Moreover, by using that $X$ and $\epsilon$ are uncorrelated it follows $\sigma_Y^2= \Var(\skalarV{X,\sol})+\sigma^2\Var(\epsilon)\leq \Ex \skalarV{X,\sol}^2+\sigma^2\leq \normV{\sol}^2\Ex\normV{X}^2 +\sigma^2$. Hence, $\sigma_Y^2\leq \br \Ex\normV{X}^2 +\sigma^2$ because  $\bw$ is monotonically non decreasing. The result follows now by combining the last estimates with the definition of the penalty, that is,  $\pen(m)= 192\sigma_Y^2\eta \delta_{m}/n$, which completes the proof of Theorem \ref{sec:gen:th1}.\hfill $\square$
\paragraph{Technical assertions.}\hfill\\[1ex]
The following lemmas gather technical results used in the proof of Theorem \ref{sec:gen:th1}.
We begin by recalling an inequality due to  \cite{Tala96}, which can be found e.g. in 
\cite{CRT2006}. 

\begin{lem}[Talagrand's Inequality]
  \label{app:th1:tala}
  Let $T_1, \ldots, T_n$ be independent $\cT$-valued random variables and $\nu^*_n(r)
  = (1/n)\sum_{i=1}^n\big[r(T_i) - \Ex[r(T_i)] \big]$, for $r$
  belonging to a countable class $\cR$ of measurable functions. Then,
  for $\epsilon> 0$,
  \begin{align*}
    \Ex[\sup_{r\in\cR} |\nu^*_n(r)|^2 &- 2(1+2\epsilon)H^2]_+ \\ &\leq
    C\left(\frac{v}{n}\exp(-K_1\epsilon\frac{nH^2}v)) +
      \frac{h^2}{n^2C^2(\epsilon)}
      \exp(-K_2C(\epsilon)\sqrt{\epsilon}\frac{nH}{h}) \right)
  \end{align*}
  with $K_1 = 1/6$, $K_2= 1/(21\sqrt{2})$, $C(\epsilon) =
  \sqrt{1+\epsilon}-1$ and $C$ a universal constant and where
  \[\sup_{r\in\R}\sup_{t\in\cT}|r(t)| \leq h,\quad
  \Ex\left[\sup_{r\in\cR}|\nu^*_n(r)|\right]\leq H, \quad \sup_{r\in\cR}
  \frac{1}{n}\sum_{i=1}^n \Var(r(T_i))\leq v.\]
\end{lem}

\begin{lem}\label{app:th1:l1} Let $\ev$ be the eigenvalues associated to $X\in\cX_{\eta}^{4}$  and 
$\Ex |Y/\sigma_Y|^4\leq \eta$. Suppose  sequences $\delta$, $\triangle$ and $M$ satisfying Assumption 
\ref{ass:reg:ext1}. Then there exists a constant $K(\Sigma, \eta,\delta_1)$ only depending on $\Sigma,\eta$ and $\delta_1$ such that 
\begin{gather}\label{app:th1:l1:e1}
\sum_{m=1}^{M_n}\Ex\Bigl(\sup_{ t\in \cB_{ m}}  |\skalarV{ t,\tPhi_{\widehat h}}_\hw|^2-6\sigma_Y^2\eta\,\frac{\delta_{m}}n\Bigr)_+
\leq  K(\Sigma, \eta,\delta_1)\, \frac{\sigma^2_Y}{n} \quad\mbox{ for all }n\geq1.
\end{gather}
\end{lem}
\begin{proof}[\noindent\textcolor{darkred}{\sc Proof.}] Given $m\in\N$ and $t\in \cB_{ m}:=\{f\in \cS_m:\normV{f}_\hw\leq 1\}$ denote
\begin{equation*}
v_{t}(Y,X):=  Y\1_{\Omega_{Y,X}} \skalarV{t,\tPhi_{X}}_\hw= \sum_{j=1}^m \frac{\hw_j[t]_j}{\lambda_j}  Y\1_{\Omega_{Y,X}} [X]_j,\end{equation*}
then it is easily seen that $\skalarV{t,\tPhi_{\widehat h}}_\hw=(1/n)\sum_{i=1}^n\{v_{t}(Y_i,X_i) - \Ex v_{t}(Y_i,X_i)\}$. Below we show the following three bounds
\begin{gather}\label{app:th1:l1:e1:1}
\sup_{ t\in \cB_{ m}}\sup_{y\in\R ,x\in L^2[0,1]}   |v_t(y,x)| \leq  \sigma_Y n^{1/3}  \delta_m^{1/2}=:h,\\\label{app:th1:l1:e1:2}
\Ex \sup_{ t\in \cB_{m}}  |\skalarV{t,\tPhi_{\widehat h}}_\hw|^2 \leq    \sigma_Y^2 \eta \,\frac{\delta_m}{n}=:H^2,
\\\label{app:th1:l1:e1:3}
\sup_{ t\in \cB_{ m}}  \frac{1}{n} \sum_{i=1}^n \Var  (v_{t}(Y_i,X_i))  \leq  \sigma_Y^2\eta \,\triangle_m =:v.
\end{gather}
From Talagrand's inequality (Lemma \ref{app:th1:tala}) with $\epsilon=1$ we obtain by combining \eqref{app:th1:l1:e1:1}-\eqref{app:th1:l1:e1:3}
\begin{multline*}
\Ex \Bigl[\sup_{ t\in \cB_{ m}}  |\skalarV{t,\tPhi_{\widehat h}}_\hw|^2-6 H^2\Bigr]\leq C \Bigl\{\frac{v}{n} \exp\Bigl( -  \frac{nH^2}{6v} \Bigr) + \frac{h^2}{n^2} \exp\Bigl( - \frac{c\, n\, H}{h}
\Bigr)\Bigr\}\\
=C \Bigl\{\frac{\sigma_Y^2\eta \,\triangle_m}{n} \exp\Bigl( -  \frac{\delta_m}{6\triangle_m} \Bigr) + \sigma_Y^2  \frac{ n^{2/3} \delta_m}{n^2 }\exp\Bigl( - {c\, \eta}\,
n^{1/6}\Bigr)\Bigr\}
\end{multline*}
with  $c=(1-1/\sqrt{2})/21$  and some numerical constant $C>0$. By using Assumption \ref{ass:reg:ext1}, that is $\delta_m/n \leq \delta_{M_n}/n \leq \delta_1$ and  $M_n/n \leq 1$, together with  $H^2= \sigma^2_Y\eta \delta_m/n$ it follows that
\begin{align*}
\sum_{m=1}^{M_n} \Ex &\Bigl[\sup_{ t\in \cB_{m}}  |\skalarV{t,\tPhi_{\widehat h}}_\hw|^2-6\sigma^2_Y\eta \delta_m/n\Bigr]\\
&\leq C \Bigl\{\frac{\sigma^2_Y\eta}{n} \sum_{m=1}^{M_n} \triangle_m \exp\Bigl( -  \frac{\delta_m}{6\triangle_m} \Bigr) + \sigma_Y^2 \delta_1 n^{2/3} \exp\Bigl( - {c\, \eta}\,
n^{1/6}\Bigr)\Bigr\}\\
&\leq C\frac{\sigma^2_Y  }{n}  \Bigl\{\eta\;\Sigma+\delta_1 \exp\Bigl( -{c\, \eta}\,
n^{1/6} + (5/3) \log n\Bigr)\Bigr\}, 
\end{align*}
where condition \eqref{sumcondi} in Assumption \ref{ass:reg:ext1} implies the last inequality. It follows that there exists a constant $K(\Sigma,\eta,\delta_1)$ only depending on $\Sigma,\eta$ and  $\delta_1$ such that 
\begin{equation*}
\sum_{m=1}^{M_n} \Ex \Bigl[\sup_{ t\in \cB_{m}}  |\skalarV{t,\tPhi_{\widehat h}}_\hw|^2-6\sigma^2_Y\eta \delta_m/n\Bigr]
\leq \frac{\sigma^2_Y}{n}  K(\Sigma,\eta,\delta_1),\quad\mbox{for all }n\geq1,
\end{equation*}
which proves the result.

Proof of \eqref{app:th1:l1:e1:1}. From
$\sup_{ t\in \cB_{ m}}  |\skalarV{t,g}_\hw|^2= \sum_{j=1}^m \hw_j[g]_j^2$ and the definition of $\Omega_{Y,X}$ follows
\begin{multline*}
  \sup_{y\in\R ,x\in L^2[0,1], t\in \cB_{ m}} 
 |v_t(y,x)|^2 = \sup_{y\in\R ,x\in L^2[0,1]} \sum_{j=1}^m \frac{\hw_j\sigma_Y^2}{\lambda_j}\, \1_{\Omega_{y,x}}\frac{y^2}{\sigma_Y^2} \frac{[x]_j^2}{\lambda_j} 
\leq  \sigma_Y^2 n^{2/3} \sum_{j=1}^m\frac{\hw_j}{\lambda_j} 
\end{multline*}
and, hence the definition of $\delta_m$ implies \eqref{app:th1:l1:e1:1}.

Proof of \eqref{app:th1:l1:e1:2}. Since $(Y_i,X_i)$, $i=1,\dotsc,n$, form  an $n$-sample of $(Y,X)$ we have
 $$\Ex \sup_{ t\in \cB_{m}}  |\skalarV{t,\tPhi_{\widehat h}}_\hw|^2= \sum_{j=1}^{m}\frac{\hw_j}{\ev_j^2} \Var\left(\frac{1}{n}\sum\limits_{i=1}^n Y_i \1_{\Omega_{Y_i,X_i}} {[X_i]_j}\right)\leq\frac{1}{n}  \sum_{j=1}^{m}\frac{\hw_j}{\ev_j^2} \Ex \left( Y \1_{\Omega_{Y,X}} {[ X]_j}\right)^2$$ 
 and hence  from   $\Ex|Y/\sigma_Y|^4\leq \eta $ and $X\in\cX_\eta^4$  it follows  that
\begin{align*}
\Ex \sup_{ t\in \cB_{N}}  |\skalarV{t,\tPhi_{\widehat h}}_\hw|^2 
&\leq\frac{\sigma_Y^2}{n}  \sum_{j=1}^{m}\frac{\hw_j}{\lambda_j}\left( \Ex  |Y/\sigma_Y|^4 \Ex  |[X]_j/\sqrt\ev_j|^4\right)^{1/2} 
\leq  \frac{\sigma_Y^2}{n}  \eta \sum_{j=1}^{m}\frac{\hw_j}{\lambda_j}. 
\end{align*}
Thereby, the definition of $\delta_m$ implies also \eqref{app:th1:l1:e1:2}.

Proof of \eqref{app:th1:l1:e1:3}. Consider $z:=(z_j)$ with $z_j:=(\hw_j[t]_j/\sqrt{\lambda_j})/ (\sum_{j=1}^m (\hw_j^2[t]_j^2/\lambda_j))^{1/2}$ and, hence $z\in\mmS^m=\{z\in \R^m, \sum_{j=1}^m z_j^2=1\}$.   Since $(Y_i,X_i)$, $i=1,\dotsc,n$, form  an $n$-sample of $(Y,X)$ it follows that
\begin{align*}
\sup_{ t\in \cB_{ m}}  \frac{1}{n} \sum_{i=1}^n \Var  (v_{t}(Y_i,X_i)) 
\leq  \sup_{ t\in \cB_{ m}}   \Ex \Bigl( Y \1_{\Omega_{Y,X}}\sum_{j=1}^m \frac{\hw_j[t]_j}{{\ev_j}}[X]_j\Bigr)^2.
\end{align*}
 Thereby, from   $\Ex|Y/\sigma_Y|^4\leq \eta $ and $X\in\cX_\eta^4$  we conclude  that
\begin{align*}
\sup_{ t\in \cB_{ m}}  \frac{1}{n} \sum_{i=1}^n \Var  (v_{t}(Y_i,X_i)) &\leq  \sup_{ t\in \cB_{ m}}   \sigma_Y^2(\Ex |  Y/\sigma_Y|^4)^{1/2} \Bigl(\Ex \Bigl|\sum_{j=1}^m \frac{\hw_j[t]_j}{\sqrt{\ev_j}}\;\frac{[X]_j}{\sqrt\ev_j}\Bigr|^4\Bigr)^{1/2}\\ 
&\leq \sigma_Y^2\eta^{1/2}  \sup_{ t\in \cB_{ m}} \sum_{j=1}^m (\omega_j^2[t]_j^2/\ev_j) \sup_{z\in\mmS_N} \Bigl(\Ex \Bigl|\sum_{j=1}^m z_j [X]_j/\sqrt\ev_j\Bigr|^4\Bigr)^{1/2}\\ 
&\leq \sigma_Y^2\eta \sup_{ t\in \cB_{ m}} \sum_{j=1}^m (\omega_j^2[t]_j^2/\lambda_j)  
\leq \sigma_Y^2\eta  \max_{1\leq j\leq
m}\hw_j/\ev_j. 
\end{align*}
Thus  the definition of $\triangle_m$ implies now \eqref{app:th1:l1:e1:3}, which completes the proof of Lemma \ref{app:th1:l1}.
\end{proof}
 \begin{lem}\label{app:th1:l2} Let $\ev$ be the eigenvalues associated to $X\in\cX_{\xi}^{24}$  and let $\Ex |Y/\sigma_Y|^{24}\leq \xi$. Suppose  sequences $\delta$, $\triangle$ and $M$ satisfying Assumption \ref{ass:reg:ext1}. Then there exists a numerical constant $C$ such that
\begin{gather}\label{app:th1:l2:e1}
\Ex \sup_{ t\in \cB_{M_n} } |\skalarV{ t,\tPhi_{\widehat f}}_\hw|^2
\leq \sqrt{2}{\, \xi\, \sigma_Y^2\, \delta_1} / n\qquad\mbox{ and }\\\label{app:th1:l2:e2}
\Ex \sup_{ t\in \cB_{M_n} } |\skalarV{ t,\hPhi_{\hg}-\tPhi_{\hg}}_\hw|^2\leq C \frac{\xi}{n}\{ \sigma_Y^2\, \delta_1  + \normV{\sol}_\hw^2\} \{1+(\Ex\normV{X}^2)^2\}\quad\mbox{for all }n\geq1.
\end{gather}
\end{lem}
\begin{proof} Since $(Y_i,X_i)$, $i=1,\dotsc,n$, form  an $n$-sample of $(Y,X)$ it follows that
\begin{multline*}
\Ex\sup_{ t\in \cB_{M_n} } |\skalarV{ t,\tPhi_{\widehat f}}_\hw|^2
= \sum_{j=1}^{M_n} \frac{\hw_j }{\lambda_j^2} \Var\left(\frac{1}{n}\sum_{i=1}^n  Y \1_{\Omega_{Y,X}^c} {[X]_j}\right)
\leq \sum_{j=1}^{M_n} \frac{\hw_j }{\lambda_j^2n} \Ex\left( Y \1_{\Omega_{Y,X}^c} {[X]_j}\right)^2.
\end{multline*}
 Thereby, from   $\Ex|Y/\sigma_Y|^{24}\leq \xi $ and $X\in\cX_\xi^{24}$  we conclude  that
 \begin{eqnarray*}
 \Ex\sup_{ t\in \cB_{M_n} } |\skalarV{ t,\tPhi_{\widehat f}}_\hw|^2&\leq &\frac{\sigma_Y^2}{n} \sum_{j=1}^{M_n} \frac{\hw_j }{\lambda_j} \left(\Ex |Y/\sigma_Y|^8  \Ex|[X]_j/\sqrt\ev_j|^8\right)^{1/4} P(\Omega_{Y,X}^c)^{1/2}\\&\leq &
 \frac{\sigma_Y^2 \xi^{1/2}}{n}\sum_{j=1}^{M_n} \frac{\hw_j}{\ev_j}  P(\Omega_{Y,X}^c)^{1/2}\leq \sigma_Y^2 \xi^{1/2}\frac{\delta_{M_n}}{n} P(\Omega_{Y,X}^c)^{1/2}
 \end{eqnarray*}
 where the last inequality follows from the property  $\delta_m\geq\sum_{j=1}^{m} \frac{\hw_j}{\ev_j} $ for all $m\geq1$. Hence by using Assumption \ref{ass:reg:ext1}, that is $\delta_{M_n}/n\leq \delta_1$, we obtain
 \begin{equation*}
 \Ex\sup_{ t\in \cB_{M_n} } |\skalarV{ t,\tPhi_{\widehat f}}_\hw|^2\leq  \sigma_Y^2 \delta_1 \xi^{1/2} \, P(\Omega_{Y,X}^c)^{1/2}.
 \end{equation*}
The estimate \eqref{app:th1:l2:e1} follows now from  $P(\Omega_{Y,X}^c)\leq 2\xi /n^2$, which can be realized as follows. 
 Since $\Omega_{Y,X}^c=\{|Y/\sigma_Y|>n^{1/6}\} \cup \bigcup_{j=1}^{M_n} \{|[X]_j/\sqrt\ev_j|>n^{1/6}\}$ it follows  by using Markov's inequality together with  $\Ex|Y/\sigma_Y|^{24}\leq \xi $ and $X\in\cX_\xi^{24}$ that 
\begin{eqnarray*}
P(\Omega_{Y,X}^c)&\leq &P(|Y/\sigma_Y|>n^{1/6})+ \sum_{j=1}^{M_n}P(|[X]_j/\sqrt\ev_j|>n^{1/6})\\
&\leq & \frac{\Ex |Y/\sigma_Y|^{18}}{n^{3}} + \sum_{j=1}^{M_n}\frac{\Ex |[X]_j/\sqrt\ev_j|^{18}}{n^3} 
\leq \frac{\xi}{n^3}(1+M_n)
\end{eqnarray*}
Thus, under Assumption \ref{ass:reg:ext1}, that is, $M_n/n\leq 1$, we obtain $P(\Omega_{Y,X}^c)\leq 2\xi /n^2$, which completes the proof of \eqref{app:th1:l2:e1}.

Proof of \eqref{app:th1:l2:e2}. Consider the decomposition
\begin{multline}\label{pr:app:th1:l2:e2:1}
{\sup_{ t\in \cB_{M_n}}  |\skalarV{ t,\hPhi_{\hg}-\tPhi_{\hg}}_\hw|^2} =  \sum_{j=1}^{M_n} \frac{\hw_j}{\ev_j}\Bigl(\frac{\ev_j}{\hev_j}\1\{\hev_j\geq1/n\}-1\Bigr)^2 \left(\frac{1}{n}\sum_{i=1}^n   Y_i \frac{[ X_i]_j}{\sqrt\ev_j}\right)^2\\
\hfill\leq 2 \sum_{j=1}^{M_n} \frac{\hw_j}{\ev_j}\Bigl(\frac{\ev_j}{\hev_j}-1\Bigr)^2\1\{\hev_j\geq1/n\} \left(\frac{1}{n}\sum_{i=1}^n  Y_i \frac{[ X_i]_j}{\sqrt\ev_j} - \sqrt{\ev_j}[\sol]_j\right)^2\\
\hfill +2 \sum_{j=1}^{M_n} \hw_j [\sol]_j^2\Bigl(\frac{\ev_j}{\hev_j}-1\Bigr)^2\1\{\hev_j\geq1/n\}\hfill\\
\hfill+2 \sum_{j=1}^{M_n} \frac{\hw_j}{\ev_j} \left(\frac{1}{n}\sum_{i=1}^n   Y_i \frac{[ X_i]_j}{\sqrt\ev_j} - \sqrt{\ev_j}[\sol]_j\right)^2\1\{\hev_j<1/n\}\\\hfill +2  \sum_{j=1}^{M_n}\hw_j [\sol]_j^2\1\{\hev_j<1/n\}
\end{multline}
where we bound  each summand separately. First, from \eqref{app:th1:l4:e1} and \eqref{app:th1:l4:e4}  in Lemma  \ref{app:th1:l4} together with $X\in\cX_{\xi}^{24}$  and $\Ex |Y/\sigma_Y|^{24}\leq \xi$ it follows that there exists a numeric constant $C>0$ such that 
\begin{multline}\label{pr:app:th1:l2:e2:2:1}
\Ex \sum_{j=1}^{M_n} \frac{\hw_j}{\ev_j}\Bigl(\frac{\ev_j}{\hev_j}-1\Bigr)^2\1\{\hev_j\geq1/n\} \left(\frac{1}{n}\sum_{i=1}^n   Y_i \frac{[ X_i]_j}{\sqrt\ev_j} - \sqrt{\ev_j}[\sol]_j\right)^2\\
\hfill\leq \sum_{j=1}^{M_n} \frac{\hw_j}{\ev_j} \Bigl[\Ex |{\ev_j}/{\hev_j}-1|^4\1\{\hev_j\geq1/n\} \Bigr]^{1/2} \Bigl[\Ex \left(\frac{1}{n}\sum_{i=1}^n   Y_i \frac{[ X_i]_j}{\sqrt\ev_j}  - \sqrt{\ev_j}[\sol]_j\right)^4 \Bigr]^{1/2}\\
\leq C \frac{\sigma_Y^2\xi}{n} \sum_{j=1}^{M_n}\frac{\hw_j}{n\ev_j} \{ \ev_j^2 +1\};
\end{multline}
\hfill\\[-8ex]
\begin{multline}\label{pr:app:th1:l2:e2:2:2}
\Ex \sum_{j=1}^{M_n} \hw_j [\sol]_j^2\Bigl(\frac{\ev_j}{\hev_j}-1\Bigr)^2\1\{\hev_j\geq1/n\}
\leq C\frac{\xi}{n}  \sum_{j=1}^{M_n}\hw_j[\sol]_j^2\{\ev_j^2+1\}.\hfill
\end{multline}
Furthermore, Assumption \ref{ass:reg:ext1} (ii), i.e., $2/n\leq \min\{\lambda_{j}:{1\leq j\leq M_n}\}$, implies $P( \hev_j<1/n)\leq P(\hev_j/\ev_j< 1/2)$. Thereby, from \eqref{app:th1:l4:e1} and \eqref{app:th1:l4:e3}  in Lemma  \ref{app:th1:l4} together with $X\in\cX_{\xi}^{24}$  and $\Ex |Y/\sigma_Y|^{24}\leq \xi$ it follows that there exists a numeric constant $C>0$ such that 
\begin{multline}\label{pr:app:th1:l2:e2:2:3}
\Ex \sum_{j=1}^{M_n} \frac{\hw_j}{\ev_j} \left(\frac{1}{n}\sum_{i=1}^n  Y_i \frac{[ X_i]_j}{\sqrt\ev_j}- \sqrt{\ev_j}[\sol]_j\right)^2\1\{\hev_j<1/n\}\\
\hfill\leq \sum_{j=1}^{M_n} \frac{\hw_j}{\ev_j}  \Bigl[\Ex \left(\frac{1}{n}\sum_{i=1}^n  Y_i\frac{[ X_i]_j}{\sqrt\ev_j}- \sqrt{\ev_j}[\sol]_j\right)^4 \Bigr]^{1/2} P( \hev_j/\ev_j<1/2)^{1/2}\\
\leq C\frac{\sigma_Y^2\xi}{n} \sum_{j=1}^{M_n}  \frac{\hw_j}{n\lambda_j};
\end{multline}
\hfill\\[-8ex]
\begin{multline}\label{pr:app:th1:l2:e2:2:4}
\Ex \sum_{j=1}^{M_n} \hw_j [\sol]_j^2\1\{\hev_j<1/n\}\leq  \sum_{j=1}^{M_n} \hw_j[\sol]_j^2 P(\hev_j/\ev_j<1/2)
\leq C \frac{\xi}{n}\sum_{j=1}^{M_n}\hw_j [\beta]_j^2.\hfill
\end{multline}
Combining the decomposition \eqref{pr:app:th1:l2:e2:1} and the bounds \eqref{pr:app:th1:l2:e2:2:1} - \eqref{pr:app:th1:l2:e2:2:4} we obtain
\begin{multline*}
\Ex\sup_{ t\in \cB_{M_n}}  |\skalarV{ t,\hPhi_{\hg}-\tPhi_{\hg}}_\hw|^2\leq C \frac{\xi}{n}\Bigl\{ \sum_{j=1}^{M_n}  \frac{\hw_j}{n\ev_j} \sigma_Y^2\{\ev_j^2+2\} + \sum_{j=1}^{M_n}\hw_j [\sol]_j^2\{\ev_j^2+2\}\Bigr\}.\hfill
\end{multline*}
Therefore the properties $\Ex\normV{X}^2\geq\max_{j\geq 1}\ev_j$ and $\delta_m\geq\sum_{j=1}^{m} \frac{\hw_j}{\ev_j} $ for all $m\geq1$ imply
\begin{multline*}
\Ex\sup_{ t\in \cB_{M_n}}  |\skalarV{ t,\hPhi_{\hg}-\tPhi_{\hg}}_\hw|^2\leq C \frac{\xi}{n}\{ \sigma_Y^2\delta_{M_n}/n  + \normV{\sol}_\hw^2\} \{(\Ex\normV{X}^2)^2+2\} .\hfill
\end{multline*}
Thus \eqref{app:th1:l2:e2} follows now  from $\delta_{M_n}/n\leq \delta_1$ (Assumption \ref{ass:reg:ext1}), which completes the proof.\end{proof}

\begin{lem}\label{app:th1:l4} Suppose $X\in\cX_{\eta_{4k}}^{4k}$ and $\Ex |Y/\sigma_Y|^{4k}\leq \eta_{4k}$, $k\geq1$. Then
for some numeric constant $C_{k}>0$ only depending on $k$ we have
\begin{align}\label{app:th1:l4:e1}
\Ex\left(\frac{1}{n}\sum_{i=1}^n    Y_i \frac{[X_i]_j}{\sqrt\ev_j} - \sqrt{\ev_j}[\sol]_j\right)^{2k}&\leq 
 C_{k} \sigma_Y^{2k} \eta_{4k}\, n^{-k},\\\label{app:th1:l4:e2}
\Ex| {\hev_j}/{\ev_j} -1|^{2k}&\leq  C_{k}\eta_{4k}\, n^{-k}.
\end{align}
If in addition $w_1\geq 2$ and $w_2\leq 1/2$, then we obtain
\begin{gather}\label{app:th1:l4:e3}
\sup_{j\in\N}P(\widehat{\lambda}_j/\lambda_j\geq w_1)\leq C_{k}\eta_{4k} \, n^{-k}\;\text{ and }\;
\sup_{j\in\N}P(\widehat{\lambda}_j/\lambda_j< w_2)\leq C_{k}  \eta_{4k}\, n^{-k}.
\end{gather}
Moreover, if $X\in\cX_{\eta_{12k}}^{12k}$, $k\geq1$, then
for some numeric constant $C_{k}>0$ only depending on $k$ we have
\begin{gather}\label{app:th1:l4:e4}
\Ex| {\ev_j}/{\hev_j} -1|^{2k}\1\{\hev_j\geq 1/n\}\leq  C_{k} \eta_{12k}\{ \ev_j^{2k} +1\}n^{-k} .
\end{gather}
\end{lem}
\begin{proof}Since $ \Ex Y [X]_j = {\ev_j}[\sol]_j$ the independence within the sample of $(Y,X)$  implies by using    Theorem 2.10 in \cite{Petrov1995}  for some generic constant $C_{k}$ that
\begin{multline*}
\Ex\left(\frac{1}{n}\sum_{i=1}^n   Y_i \frac{[X_i]_j}{\sqrt\ev_j} - \sqrt{\ev_j}[\sol]_j\right)^{2k}\leq C_{2k}\sigma_Y^{2k} n^{-k} 
\Ex |Y/\sigma|^{2k} |[X]_j/\sqrt\ev_j|^{2k} \\\leq  C_{k} \sigma_Y^{2k}n^{-k} \left(\Ex |Y/\sigma|^{4k} \Ex|[X]_j/\sqrt\ev_j|^{4k}\right)^{1/2}.
\end{multline*}
Then the last estimate together with $X\in\cX_{\eta_{4k}}^{4k}$ and $\Ex |Y/\sigma_Y|^{4k}\leq \eta_{4k}$ implies \eqref{app:th1:l4:e1}. Furthermore,  since $\{(|[X_{i}]_j|^2/\ev_j-1)_i\}$  are independent and identically distributed with mean zero, it follows by applying again Theorem 2.10 in \cite{Petrov1995} that $\Ex| {\hev_j}/{\ev_j} -1|^{2k}\leq  C_{k} n^{-k} \Ex||[X]_j/\sqrt\ev_j|^{2}-1|^{2k}$. Thus, the condition  $X\in\cX_{\eta_{4k}}^{4k}$ implies \eqref{app:th1:l4:e2}. 

Proof of \eqref{app:th1:l4:e3}. If $w\geq 2$ then  $P(\hev_j/\ev_j\geq w)\leq P(|\hev_j/\ev_j-1|\geq 1)$. Thus  applying Markov's inequality together with \eqref{app:th1:l4:e2} implies the first bound in \eqref{app:th1:l4:e3}, while the second follows in analogy.

Proof of \eqref{app:th1:l4:e4}. By using twice the elementary inequality $| {\hev_j}/{\ev_j} -1|^{2k} + | {\hev_j}/{\ev_j}|^{2k}\geq 1/2^{2k-1} $ we conclude that 
\begin{multline*}
\Ex| {\ev_j}/{\hev_j} -1|^{2k}\1\{\hev_j\geq 1/n\}\leq 2^{2k-1}\{\Ex | {\hev_j}/{\ev_j} -1|^{4k}\frac{\ev_j^{2k}}{\hev_j^{2k}}\1\{\hev_j\geq 1/n\}  +  \Ex| {\hev_j}/{\ev_j} -1|^{2k}\}\\
\hfill \leq 2^{4k-2} \ev_j^{2k} n^{2k} \Ex | {\hev_j}/{\ev_j} -1|^{6k} +   2^{4k-2} \Ex| {\hev_j}/{\ev_j} -1|^{4k}  + 2^{2k-1} \Ex| {\hev_j}/{\ev_j} -1|^{2k}\}.
\end{multline*}
Thus, \eqref{app:th1:l4:e4} follows from \eqref{app:th1:l4:e2} since  $X\in\cX_{\eta_{12k}}^{12k}$,  which proves the lemma.\end{proof}

\subsection{Proof of Proposition \ref{ratesknownlambda1}}
\paragraph{Case [P-P]}
Since $2a+2s+1> 0$ it follows that the sequences $\delta, \Delta$ and $M$ with  $\delta_m\asymp m^{2a+2s+1}$,
$\Delta_m\asymp m^{(2a+2s)_{\vee0}}$ and $M_n\asymp
n^{1/(2a+1+(2s)_{\vee 0})}$, respectively, satisfy 
Assumption \ref{ass:reg:ext1}. Note that $\delta_{M_n}/n\leq1$, $M_n/n\leq 1$, $\min_{1\leq j\leq  M_n} \ev_{j} \geq 2/n$  and $\forall C>0$,
$$\sum_m \triangle_m\exp(-C\delta_m/\Delta_m)
\asymp \sum_m m^{(2a+2s)_{\vee0}}\exp(-Cm^{(2a+2s+1)\wedge 1}) <+\infty.$$
Therefore we can apply Theorem \ref{sec:gen:th1} and hence Corollary \ref{sec:gen:prop1:cor}. In particular, by using  $\mdagn\asymp n^{1/(2a+2p+1)}$, which satisfies $\bw_{\mdagn}\delta_{\mdagn}/(n\hw_{\mdagn}) \asymp 1$,  it follows that the adaptive
estimator $\widehat\beta_{\widehat m}$ reaches the optimal rate $\hw_{\mdagn}/{\bw_{\mdagn}}
\asymp n^{-2(p-s)/(2p+2a+1)}$.

\paragraph{Case [E-P]} The sequences $\delta, \Delta, M$ are unchanged w.r.t. the previous case [P-P] and  hence Assumption \ref{ass:reg:ext1} is still satisfied. From Corollary \ref{sec:gen:prop1:cor} follows now again that the adaptive
estimator $\widehat\beta_{\widehat m}$ attains the optimal rate $\hw_{\mdagn}/{\bw_{\mdagn}}\asymp n^{-1} (\log
n)^{(2a+1+2s)/(2p)}$  since $\mdagn\asymp \{\log [n (\log n)^{-(2a+1)/(2p)}] \}^{1/(2p)}$  satisfies
$\bw_{\mdagn}\delta_{\mdagn}/(n\hw_{\mdagn}) \asymp 1$. 

\paragraph{Case [P-E]} Consider the sequences $\delta, \Delta$ and $M$ with  $\delta_m=m^{2a+1+(2s)_{\vee0}}\exp(m^{2a})$, $\Delta_m=m^{(2s)_{\vee0}}\exp(m^{2a})$ and $M_n=(\log
n/(\log n)^{2a+1+(2s)_{\vee0})/(2a)})^{1/(2a)}$ respectively. Then Assumption
\ref{ass:reg:ext1} is satisfied, that is  $\delta_{M_n}/n\leq 1$,
$M_n/n\leq 1$, $\min_{1\leq j\leq  M_n} \ev_{j} \geq 2/n$  and $\forall C>0$,
$$\sum_m \Delta_m\exp(-C\delta_m/\Delta_m) \leq \sum_m  m^{(2s)_{\vee0}}\exp(m^{2a}) \exp(-C m^{2a+1}) <+\infty.$$
Moreover, $\bw_{\mdag}\delta_{\mdag}/(n\hw_{\mdag})\asymp 1$ implies $\mdagn\asymp (\log n/(\log n)^{(2a+2p+1)/(2a)})^{1/(2a)}$.  Finally, due to  Corollary \ref{sec:gen:prop1:cor}  the adaptive
estimator $\widehat\beta_{\widehat m}$ attains again the optimal rate ${\hw_{\mdag}}/{\bw_{\mdag}}\asymp(\log
n)^{-(p-s)/a}$, which completes the proof of Proposition \ref{ratesknownlambda1}.\hfill$\square$

\subsection{Proof of Proposition \ref{ratesknownlambda2}}
Let $\Delta_m:=\max_{1\leq j\leq m}\hw_j/\ev_j$, $\kappa_m:=\max_{1\leq j\leq m}(\hw_j)_{\vee1}/\ev_j$ and
$\delta_m:= m
\Delta_m\Bigl|\frac{\log (\kappa_m\vee (m+2))}{\log(m+2)}\Bigr|$ as defined in (\ref{deltamref}). Note that $|\log
(\kappa_m\vee (m+2))/\log(m+2)|\geq 1$ and hence $\delta_m\geq
\sum_{j=1}^m\hw_j/\ev_j.$

\paragraph{Case [P-P] and [E-P].} Since $a+s\geq 0$ it is easily verified that $\Delta_m\asymp m^{2a+2s}$, $\kappa_m\asymp m^{2a+(2s)_{\vee0}}$ with $ |\log (\kappa_m\vee (m+2))/\log(m+2)|\asymp (2a+(2s)_{\vee0})>1$  and hence,
$\delta_m \asymp  m^{1+2a+2s}$. Therefore, the result follows from Proposition \ref{ratesknownlambda1} case [P-P] and [E-P] since  both sequences $\delta$ and $\Delta$ are unchanged.

\paragraph{Case [P-E]} We have  $\Delta_m\asymp m^{2s}\exp(m^{2a})$,
$\kappa_m\asymp m^{(2s)_{\vee0}}\exp(m^{2a})$ with, for all $m$ sufficiently large, 
$\log( \kappa_m \vee (m+2))/\log(m+2)| \asymp m^{2a}\frac{(1+(2s)_{\vee0}(\log m)m^{-2a})}{\log( m+2) }$ and hence $ 
\delta_m \asymp m^{1+2a+2s}\exp(m^{2a}) \frac{(1+(2s)_{\vee0}(\log
m)m^{-2a})}{\log m }.$ Then straightforward calculus shows that  Assumption \ref{ass:reg:ext1} (i) is fulfilled. 
Moreover, consider the sequence $M$ given in Assumption \ref{ass:reg:ext1} (ii), where 
$M_n\asymp (\log  \frac{n\;(\log\log n)/(2a)}{(\log
n)^{(1+2a+(2s)_{\vee0})/(2a)}})^{1/(2a)} = (\log
n)^{1/(2a)}\Bigl(1+o(1)\Bigr)
$, then also Assumption \ref{ass:reg:ext1} (ii) is satisfied (as in the proof of case [P-E] in Proposition \ref{ratesknownlambda1}). Due to  Corollary \ref{sec:gen:prop1:cor} it remains to balance $n\asymp\bw_{\mdag}\delta_{\mdag}/\hw_{\mdag}\asymp
(\mdag)^{1+2a+2p}\exp((\mdag)^{2a})/(\log \mdag)$ which implies
$\mdagn\asymp (\log \frac{n\;(\log\log n)/(2a)}{(\log n)^{(1+2a+2p)/(2a)}})^{1/(2a)} = (\log n)^{1/(2a)}\Bigl(1+o(1)\Bigr)$. Hence,  $ \hw_{\mdag}/\bw_{\mdag}\asymp  (\log n)^{(p-s)/a} $ is the rate attained by  the adaptive
estimator $\widehat\beta_{\widehat m}$ which is  optimal and completes the proof of Proposition \ref{ratesknownlambda2}.\hfill$\square$

\subsection{Proof of Theorem \ref{sec:gen:th2}}
We begin by defining additional notations to be used in the proof.
Consider sequences  $\delta$, $\triangle$, $M$ and $\mdag$ satisfying  Assumption \ref{ass:reg:ext2} and  the random upper bound $\hM$  defined in \eqref{sec:gen:rub}. Denote by
$\Omega:=\Omega_I\cap\Omega_{II}$ the event given by
\begin{align*}
\Omega_I&:=\left\{ \forall j\in \{1, \dots, M_n\}, \left|\frac 1{\hev_j}-\frac 1{\ev_j}\right| < \frac 1{2\ev_j} \mbox{ and } \hev_j\geq 1/n\right\},\\
\Omega_{II}&:=\{ \mdagn \leq \hM_n\leq M_n\}.
\end{align*}
It is easily seen that on $\Omega_I$ we have  for all   $1\leq m\leq
M_n$ $$ (1/2)\Delta_m\leq
\widehat\Delta_m\leq(3/2)\Delta_m\qquad\mbox{and}\qquad(1/2)\kappa_m\leq
\widehat\kappa_m\leq(3/2)\kappa_m$$ and hence $(1/2)[\kappa_m\vee(m+2)]\leq[
\widehat\kappa_m\vee (m+2)]\leq(3/2)[\kappa_m\vee(m+2)]$ which implies
\begin{multline*}(1/2)  m  \Delta_m\Bigl(\frac{ \log[ \kappa_m\vee (m+2)]}{\log(m+2)}\Bigr)\Bigl(1-\frac{\log 2}{\log (m+2)}\frac{\log (m+2)}{\log (\kappa_m \vee [m+2])}\Bigr) \leq
\widehat\delta_m  \\\hfill\leq (3/2)   m  \Delta_m\Bigl(\frac{\log (\kappa_m \vee [m+2])}{\log( m+2)}\Bigr) \Bigl(1+\frac{\log3/2}{\log(m+2)
} \frac{\log (m+2)}{\log (\kappa_m \vee [m+2])}\Bigr),\end{multline*}
 together with $ {\log (\kappa_m \vee [m+2])}/{\log( m+2)}\geq 1$ we get
 \begin{multline*} \delta_m/10\leq (\log 3/2)/(2 \log 3) \delta_m\leq (1/2)   \delta_m[1-(\log 2)/\log(m+2)] \leq \widehat\delta_m \\ \leq (3/2)  \delta_m [1+  (\log3/2)/\log (m+2)]\leq 3\delta_m.\end{multline*}
 Since $\pen(m)= 192 \sigma_Y^2\eta \delta_m n^{-1}$ and $\hpen(m)= 1920 \sigma_Y^2\eta \widehat\delta_m n^{-1}$ 
 it follows on $\Omega_I$ that  $\pen(m)\leq \hpen(m)\leq 30 \pen(m)$ for all   $1\leq m\leq
M_n$, and hence
\begin{multline*}
\Bigl( \pen(\mdagn\vee \whm) + \hpen(\mdagn)-\hpen(\whm) \Bigr)\1_{\Omega} \leq \Bigl( \pen(\mdagn)+\pen(\whm) + \hpen(\mdagn)-\hpen(\whm) \Bigr)\1_{\Omega}\\
\hfill\leq 31\pen(\mdagn)
\end{multline*}
 by using   $1\leq \whm\leq \hM_n$ and  $\mdagn\leq \hM_n\leq M_n$. On the other hand,  it is not hard to see that on $\Omega_I^c$ we have 
 $\widehat\Delta_m\leq n \max_{1\leq j\leq m}\hw_j$ and $\widehat\kappa_m\leq n$ for all   $m\geq
1$. From these properties we conclude that for all $1\leq m\leq M_n$
 \begin{equation}\label{app:proofs:gen:2:def:3} \widehat\delta_m  \leq  m n (\max_{1\leq j\leq m}\hw_j) \frac{\log (n \vee(m+2))}{\log(m+2)}\leq m n (\max_{1\leq j\leq m}\hw_j)\log(n+2),
 \end{equation}
 which implies  $\hpen(\mdagn)\leq 1920 \sigma_Y^2\eta M_n (\max_{1\leq j\leq M_n}\hw_j) {\log (n +2)} $ and hence
\begin{multline}\label{app:proofs:gen:2:def:4}
\Bigl( \pen(\mdagn\vee \whm) + \hpen(\mdagn)-\hpen(\whm) \Bigr)\1_{\Omega_I^c\cap\Omega_{II}} \\\hfill\leq \Bigl( \pen(M_n) + 1920 \sigma_Y^2\eta M_n (\max_{1\leq j\leq M_n}\hw_j) {\log (n+2)}\Bigr)\1_{\Omega_I^c\cap\Omega_{II}}\hfill\\
\hfill \leq 1920 \sigma_Y^2\eta \Bigl( \delta_{M_n}/n + M_n (\max_{1\leq j\leq M_n}\hw_j) {\log (n+2)}\Bigr)\1_{\Omega_I^c\cap\Omega_{II}}
\end{multline}
We shall prove in the end of this section the technical Lemma \ref{app:th2:l1}  which is used in the following steps of the proof together with the technical Lemmas \ref{app:th1:l1} - \ref{app:th1:l4} above.\\[1ex]

Consider now the decomposition
\begin{multline}\label{pr:th2:e1}
\Ex\normV{\hsol_{\whm}-\sol}_\hw^2= \Ex\normV{\hsol_{\whm}-\sol}_\hw^2\1_\Omega + \Ex\normV{\hsol_{\whm}-\sol}_\hw^2\1_{\Omega_I^c\cap\Omega_{II}} +  \Ex\normV{\hsol_{\whm}-\sol}_\hw^2\1_{\Omega_{II}^c}.
\end{multline}
Below we show that there exist a numerical constant $C'>0$ and a constant $K'=K'(\Sigma, \eta,\xi,\delta_1)$  only depending on $\Sigma, \eta,\xi$ and $\delta_1$  such that for all $n\geq1$ we have
\begin{gather}\label{pr:th2:e1:1}
\Ex \normV{\hsol_{\whm}-\sol}^2_\hw\1_\Omega\leq C'\Bigl\{ \normV{\sol-\sol_{\mdagn}}^2_\hw + \frac{\delta_{\mdagn}}{n} \,\sigma_Y^2\eta  
+ \frac{K'}{n}\,  \sigma_Y^2\, [ \delta_1 + \normV{\sol}^2_\hw][ 1 +
(\Ex\normV{X}^2)^2] \Bigr\}\\\label{pr:th2:e1:2}
\Ex \normV{\hsol_{\whm}-\sol}^2_\hw\1_{\Omega_I^c\cap\Omega_{II}}\leq C'\Bigl\{ \normV{\sol-\sol_{\mdagn}}^2_\hw +  \frac{K'}{n}\,  \sigma_Y^2\, [ \delta_1+ \normV{\sol}^2_\hw][ 1 +
(\Ex\normV{X}^2)^2] \Bigr\}\\\label{pr:th2:e1:3}
\Ex\normV{\hsol_{\whm}-\sol}^2_\hw\1_{\Omega_{II}^c}      \leq C' \frac{\xi}{n} \, [ \sigma_Y^2 + \normV{\sol}^2_\hw][ 1 +
\Ex\normV{X}^2].
\end{gather}
Since  $(\hw/\bw)$ is monotonically non increasing we obtain in case $\sol\in\cF_\bw^\br$ that  $\normV{\sol}^2_\hw\leq \br$ and $\normV{\sol-\sol_{\mdagn}}^2_\hw\leq (\hw_{\mdagn}/\bw_{\mdagn})\br$. Moreover, we have  $\sigma_Y^2\leq \br \Ex\normV{X}^2 +\sigma^2$. From these properties by combining the decomposition \eqref{pr:th2:e1} and the estimates \eqref{pr:th2:e1:1} - \eqref{pr:th2:e1:3} we conclude that there exists a numerical constant $C>0$ and a constant $K=K(\Sigma, \eta,\xi,\delta_1)$  only depending on $\Sigma, \eta,\xi$ and $\delta_1$  such that for all $n\geq1$ 
\begin{equation*}\Ex\normV{\hsol_{\whm}-\sol}_\hw^2\leq C\Bigl\{ \frac{\hw_{\mdagn}}{\bw_{\mdagn}}\br + \frac{\delta_{\mdagn}}{n} \,[\br \Ex\normV{X}^2 +\sigma^2]\eta  + \frac{K}{n}\,  [\br \Ex\normV{X}^2 +\sigma^2]\, [ 1+\delta_1 + \br][ 1 + 
(\Ex\normV{X}^2)^2] \Bigr\} .\end{equation*}
The result follows now  from the definition of $\mdagn$, that is,  ${\bw_{\mdagn}\delta_{\mdagn}}/({n\,\hw_{\mdagn}})\leq \uc$.

Proof of \eqref{pr:th2:e1:1}. Observe that on $\Omega$ we have
$\mdag_n\leq \hM_n\leq M_n$. Thus, following line by line the proof
of \eqref{pr:th1:e3}  it is easily seen that
\begin{eqnarray*}\nonumber
(1/4) \normV{\hsol_{\whm}-\sol}^2_\hw\1_\Omega &\leq & (7/4) \normV{\sol-\sol_{\mdagn}}^2_\hw + 32\sum_{m=1}^{M_n}\Bigl(\sup_{ t\in \cB_{ m}}  |\skalarV{ t,\tPhi_{\widehat h}}_\hw|^2-6\sigma_Y^2\eta \delta_m/n\Bigr)_+\\ \nonumber && +   32 \sup_{ t\in \cB_{M_n} } |\skalarV{ t,\tPhi_{\widehat f}}_\hw|^2+ { 32
\sup_{ t\in \cB_{M_n} } |\skalarV{ t,\hPhi_{\hg}-\tPhi_{\hg}}_\hw|^2} \\ \nonumber && +    \Bigl( \pen(\mdagn\vee \whm) + \hpen(\mdagn)-\hpen(\whm) \Bigr)\1_{\Omega}\\ \nonumber
&\leq & (7/4) \normV{\sol-\sol_{\mdagn}}^2_\hw + 32\sum_{m=1}^{M_n}\Bigl(\sup_{ t\in \cB_{ m'}}  |\skalarV{ t,\tPhi_{\widehat h}}_\hw|^2-6\sigma_Y^2\eta \delta_m/n\Bigr)_+\\  \nonumber &&
+   32 \sup_{ t\in \cB_{M_n} } |\skalarV{ t,\tPhi_{\widehat f}}_\hw|^2+ { 32
\sup_{ t\in \cB_{M_n} } |\skalarV{ t,\hPhi_{\hg}-\tPhi_{\hg}}_\hw|^2} \\ 
&&+ 4 \pen(\mdagn),
\end{eqnarray*}
where the last inequality follows from \eqref{app:proofs:gen:2:def:3}. Combining the last bound with \eqref{app:th1:l1:e1}  in Lemma \ref{app:th1:l1},  \eqref{app:th1:l2:e1}  and \eqref{app:th1:l2:e2}  in Lemma \ref{app:th1:l2} we conclude  that there exists a numerical constant $C'>0$ and a constant $K'=K'(\Sigma, \eta,\xi,\delta_1)$   depending on $\Sigma, \eta,\xi,\delta_1$  only such that \eqref{pr:th2:e1:1} for all $n\geq1$ holds true.

Proof of \eqref{pr:th2:e1:2}. Note that on $\Omega_I^c\cap\Omega_{II}$ we have still $\mdag_n\leq \hM_n\leq M_n$. Thus, by using   \eqref{app:proofs:gen:2:def:4} rather than \eqref{app:proofs:gen:2:def:3} it follows in analogy to \eqref{pr:th2:e1} that
\begin{multline*}
(1/4) \normV{\hsol_{\whm}-\sol}_\hw^2\1_{\Omega_I^c\cap\Omega_{II}}\leq (7/4) \normV{\sol-\sol_{\mdagn}}^2_\hw + 32\sum_{m=1}^{M_n}\Bigl(\sup_{ t\in \cB_{ m}}  |\skalarV{ t,\tPhi_{\widehat h}}_\hw|^2-6\sigma_Y^2\eta \delta_m/n\Bigr)_+\hfill\\
\hfill+   32 \sup_{ t\in \cB_{M_n} } |\skalarV{ t,\tPhi_{\widehat f}}_\hw|^2+ { 32
\sup_{ t\in \cB_{M_n} } |\skalarV{ t,\hPhi_{\hg}-\tPhi_{\hg}}_\hw|^2} +    \Bigl( \pen(\mdagn\vee \whm) + \hpen(\mdagn)-\hpen(\whm) \Bigr)\1_{\Omega_I^c\cap\Omega_{II}}\\
\hfill\leq (7/4) \normV{\sol-\sol_{\mdagn}}^2_\hw + 32\sum_{m=1}^{M_n}\Bigl(\sup_{ t\in \cB_{ m}}  |\skalarV{ t,\tPhi_{\widehat h}}_\hw|^2-6\sigma_Y^2\eta \delta_m/n\Bigr)_+\hfill\\
\hfill+   32 \sup_{ t\in \cB_{M_n} } |\skalarV{ t,\tPhi_{\widehat f}}_\hw|^2+ { 32
\sup_{ t\in \cB_{M_n} } |\skalarV{ t,\hPhi_{\hg}-\tPhi_{\hg}}_\hw|^2}\\ + 1920 \sigma_Y^2\eta \Bigl( \delta_{M_n}/n + M_n (\max_{1\leq j\leq M_n}\hw_j) {\log (n+2)}\Bigr)\1_{\Omega_I^c\cap\Omega_{II}}.
\end{multline*}
From the last bound together with \eqref{app:th1:l1:e1}  in Lemma \ref{app:th1:l1},  \eqref{app:th1:l2:e1}  and \eqref{app:th1:l2:e2}  in Lemma \ref{app:th1:l2} we conclude that there exist a numerical constant $C>0$ and a constant $K=K(\Sigma, \eta,\xi,\delta_1)$   depending on $\Sigma, \eta,\xi$ and $\delta_1$  only such that for all $n\geq1$ we have
\begin{multline}\label{pr:th2:e3}
\Ex \normV{\hsol_{\whm}-\sol}^2_\hw\1_{\Omega_I^c\cap\Omega_{II}}\leq C\Bigl\{ \normV{\sol-\sol_{\mdagn}}^2_\hw +  \frac{K}{n}\,  \sigma_Y^2\, [ \delta_1 + \normV{\sol}^2_\hw][ 1 +
(\Ex\normV{X}^2)^2] \\\hfill +\sigma_Y^2\eta \Bigl( n^{-1}\delta_{M_n} + n^{-2} M_n (\max_{1\leq j\leq M_n}\hw_j) \Bigr)n^2{\log (n+2)}P(\Omega_I^c\cap\Omega_{II})\Bigr\}.
\end{multline}
Since $X\in\cX^{24}_{\xi}$ and  $\Omega_I^c\cap\Omega_{II}\subset \{\exists  j\in\{1,\dotsc,M_n\}: |\ev_j/\hev_j -1 |> 1/2\mbox{ or } \hev_j<1/n\}$ it follows from \eqref{app:th2:l1:e1} in Lemma \ref{app:th2:l1:e1} that $P(\Omega_I^c\cap\Omega_{II})\leq C \xi\,M_n\, n^{-6}$ for some numerical constant $C>0$. Moreover, due to Assumption \ref{ass:reg:ext2} we have $\delta_{M_n}/n\leq \delta_1 $, $M_n/n\leq 1$ and $\max_{1\leq j\leq M_n}\hw_j\leq \max_{1\leq j\leq N_n}\hw_j\leq n$. Combining the last estimates and \eqref{pr:th2:e3} implies now \eqref{pr:th2:e1:2}.

Proof of \eqref{pr:th2:e1:3}. Let $\breve\sol_{m}:=\sum_{j=1}^m [\sol]_j \1\{\hev_j\geq 1/n\}\ef_j$. Then it is not hard to see that $\normV{\hsol_{m}-\breve\sol_{m}}_\hw^2 \leq \normV{\hsol_{m'}-\breve\sol_{m'}}_\hw^2$ for all $m\leq m'$ and $\normV{\breve\sol_{m}-\sol}_\hw^2\leq \normV{\sol}_\hw^2$. By using  these properties together with $1\leq \whm\leq \hM_n\leq N_n$ we conclude
\begin{eqnarray*}
\Ex\normV{\hsol_{\whm}-\sol}_\hw^2\1_{\Omega_{II}^c} &\leq & 2\{\Ex\normV{\hsol_{\whm}-\breve\sol_{\whm}}_\hw^2\1_{\Omega_{II}^c} + \Ex\normV{\breve\sol_{\whm}-\sol}_\hw^2\1_{\Omega_{II}^c}\}\\
&\leq & 2\{ \Ex \normV{\hsol_{N_n}-\breve\sol_{N_n}}_\hw^2\1_{\Omega_{II}^c} + \normV{\sol}_\hw^2 P(\Omega_{II}^c)\}.
\end{eqnarray*}
Since $X\in\cX^{28}_{\xi}$ and  $\Omega_{II}^c= \{\hM_n< \mdagn\}\cup \{\hM_n> M_n\}$ it follows from
 \eqref{app:th2:l1:e2} and \eqref{app:th2:l1:e3} in Lemma \ref{app:th2:l1} that $P(\Omega_{II}^c)\leq C \xi n^{-6}$ for some numerical constant $C>0$ and hence
 \begin{multline}\label{pr:th2:e4}
\Ex\normV{\hsol_{\whm}-\sol}_\hw^2\1_{\Omega_{II}^c}
\leq 2\{ \Ex \normV{\hsol_{N_n}-\breve\sol_{N_n}}_\hw^2\1_{\Omega_{II}^c} + C\xi\, \normV{\sol}_\hw^2 n^{-6}\}.\hfill
\end{multline}
 Moreover, from \eqref{app:th1:l4:e1} and \eqref{app:th1:l4:e2}  in Lemma  \ref{app:th1:l4} together with $X\in\cX_{\xi}^{28}$  and $\Ex |Y/\sigma_Y|^{28}\leq \xi$ it follows that there exists a numerical constant $C>0$ such that
\begin{multline}\label{pr:th2:e5}
\Ex \normV{\hsol_{N_n}-\breve\sol_{N_n}}_\hw^2\1_{\Omega_{II}^c} \leq 2n^2\sum_{j=1}^{N_n}\hw_j\Bigl\{\Ex([\hg]_j -\ev_j [\sol]_j)^2\1_{\Omega_{II}^c} + \Ex(\ev_j [\sol]_j-\hev_j [\sol]_j)^2\1_{\Omega_{II}^c}\Bigr\}\\
\hfill\leq 2  n^2  \Bigl\{  \max_{1 \leq j\leq N_n} \hw_j  \sum_{j=1}^{N_n} \ev_j \Bigl[\Ex \left(\frac{1}{n}\sum_{i=1}^n Y_i \frac{[X_i]_j}{\sqrt\ev_j} - \sqrt{\ev_j}[\sol]_j\right)^{4}\Bigr]^{1/2} P(\Omega_{II}^c)^{1/2}\\
\hfill+ \max_{j\geq 1} \ev_j\sum_{j=1}^{N_n}\hw_j  [\sol]_{j}^2[\Ex (\hev_j/\ev_j-1)^4]^{1/2}P(\Omega_{II}^c)^{1/2} \Bigr\}  \\
 \leq C \xi n^2 \Bigl\{  n^{-4} \sigma_Y^2 \max_{1 \leq j\leq N_n} \hw_j  \sum_{j\geq 1} \ev_j +  n^{-4}\max_{j\geq 1} \ev_j\normV{\sol}_\hw^2\Bigr\}.
\end{multline}
By combination of \eqref{pr:th2:e4}, \eqref{pr:th2:e5} and $\Ex\normV{X}^2=\sum_{j\geq 1}\ev_j\geq \max_{ j\geq 1}\ev_j$ we obtain
\begin{equation*}\Ex\normV{\hsol_{\whm}-\sol}^2_\hw\1_{\Omega_{II}^c}      \leq C' \Bigl\{  n^{-2} \sigma_Y^2 \xi \max_{1 \leq j\leq N_n} \hw_j  \Ex\normV{X}^2 +   \xi \{1+\Ex\normV{X}^2\} \normV{\sol}_\hw^2n^{-2}\Bigr\},\end{equation*}
for some numerical constant $C'>0$. The estimate \eqref{pr:th2:e1:3} follows now from $\max_{1 \leq j\leq N_n}\hw_j \leq  n$ (Assumption \ref{ass:reg:ext2}), which completes the proof of Theorem \ref{sec:gen:th2}.\hfill$\square$

\paragraph{Technical assertions.}\hfill\\[1ex]
The following lemma gathers technical results used in the proof of Theorem \ref{sec:gen:th2}.
\begin{lem}\label{app:th2:l1} Suppose $X\in\cX_{\eta_{4k}}^{4k}$, $k\geq 1$, with associated sequence $\ev$ of eigenvalues. Let $M$ and $\mdag$ be  sequences  satisfying  Assumption \ref{ass:reg:ext2}.  Then there exist a numerical constant $C_{k}>0$ only depending on $k$ such that for all $n\geq 1$ we have
\begin{gather}\label{app:th2:l1:e1}
P(\{\exists  j\in\{1,\dotsc,M_n\}: |\ev_j/\hev_j -1 |> 1/2\mbox{ or } \hev_j<1/n\})\leq C_{k}  \eta_{4k}\,M_n\, n^{-k},\\\label{app:th2:l1:e2}
P( \hM_n < \mdagn )\leq C_{k}  \eta_{4k}\, n^{-k}\qquad\mbox{ and}\\\label{app:th2:l1:e3}
P( \hM_n >M_n)\leq C_{k}  \eta_{4k}\, n^{-k+1}\qquad\mbox{ for all }n\geq 1.
\end{gather}
\end{lem}
\begin{proof} Proof of \eqref{app:th2:l1:e1}. We start our proof with the observation that the event $\{|\ev_j/\hev_j -1 |> 1/2\}$ can  equivalently be written as
$\{1-\hev_j/\ev_j> 1/3\mbox{ or } \hev_j/\ev_j-1>1\}$, and hence is a subset of  $ \{ |\hev_j/\ev_j -1 |> 1/3\}$. Moreover, since  $\ev_j\geq 2/n$ for all $1\leq j\leq M_n$ it follows that $\{\hev_j<1/n\}\subset \{|\hev_j/\ev_j-1|>1/2\}$. Combining both estimates we conclude
\begin{multline*}
P(\{\exists  j\in\{1,\dotsc,M_n\}: |\ev_j/\hev_j -1 |> 1/2\mbox{ or } \hev_j<1/n\})\\\hfill\leq \sum_{j=1}^{M_n}\{ P( |\hev_j/\ev_j -1 |> 1/3 ) + P(|\hev_j/\ev_j -1 |> 1/2) \} \leq 2 \sum_{j=1}^{M_n} P( |\hev_j/\ev_j -1 |> 1/3 ).
\end{multline*}
Thus  applying Markov's inequality together with
\eqref{app:th1:l4:e2} in Lemma \ref{app:th1:l4}  implies
\eqref{app:th2:l1:e1}.

Proof of \eqref{app:th2:l1:e2}. Due to the definition of $\hM_n$ given in \eqref{sec:gen:rub} the event  $\{ \hM_n < \mdagn\}$ is a subset of  $\{ \forall m\in \{\mdagn, \dotsc,n\} : \hev_m/(\hw_m)_{\vee1}< m (\log n)/n\}$ and hence $P(\hM_n < \mdagn)\leq P(\hev_{\mdagn}/\ev_{\mdagn}< 1/2)$ since $\min_{1\leq m\leq \mdagn}\ev_m/[m(\hw_m)_{\vee1}]\geq 2 (\log n)/n$ (Assumption \ref{ass:reg:ext2} (iii)). Thereby,  \eqref{app:th2:l1:e2} follows from the second bound in \eqref{app:th1:l4:e3} in Lemma \ref{app:th1:l4:e2}.

Proof of \eqref{app:th2:l1:e3}. Due to the definition \eqref{sec:gen:rub} of $\hM_n$ for $m>M_n$ the event  $\{ \hM_n = m\}$ is a subset of $\{\hev_{m}/(\hw_m)_{\vee1}\geq m (\log n)/n\}$ and hence $P(\hM_n >M_n)\leq \sum_{j=M_n+1}^{N_n}P(\hev_{m}/\ev_{m}\geq  2)$ since $2\max_{m> M_n} \ev_m/[m(\hw_m)_{\vee1}] \leq  (\log n)/n$ (Assumption \ref{ass:reg:ext2} (ii)). Thereby, the first bound in \eqref{app:th1:l4:e3} in Lemma \ref{app:th1:l4:e2} together with  $N_n/n\leq 1$  (Assumption \ref{ass:reg:ext2} (iv)) implies \eqref{app:th2:l1:e3}, which completes the proof of Lemma \ref{app:th2:l1}.\end{proof}
\subsection{Proof of Corollary \ref{optimadapt}}
First, note that in all three cases, the sequences $\delta, \Delta$, $M$ and $\mdag$ have been calculated  in the proof of Proposition \ref{ratesknownlambda2}. If in addition Assumption 
\ref{ass:reg:ext2} holds true, then from Theorem \ref{sec:gen:th2} follows that  the fully adaptive estimator attains the rate $\hw_{\mdagn}/\bw_{\mdagn}$, which  in the proof of Proposition \ref{ratesknownlambda2} has been confirmed to be optimal in all three cases. Therefore it only remains to check  $(i)$-$(iii)$ of Assumption 
\ref{ass:reg:ext2}. 
\paragraph{Case [P-P]} In this case, we have $M_n\asymp
n^{1/(2a+1+(2s)_{\vee 0})}$ and $\mdagn\asymp n^{1/(2a+2p+1)}$. 
Then $(i)$ of Assumption \ref{ass:reg:ext2} holds true, since $\min_{1\leq j\leq  M_n} \ev_{j}\asymp
M_n^{-2a}\asymp n^{-2a/(2a+1+(2s)_{\vee0})} \geq 2/n$ and 
 \begin{equation*}\max_{m\geq  M_n} \frac{\ev_{m}}{m(\hw_m)_{\vee1}}\asymp
M_n^{-1-2a-(2s)_{\vee0}}\asymp
n^{-(2a+1+(2s)_{\vee0})/(1+2a+(2s)_{\vee0})} \leq (\log n)/(2n).\end{equation*}
Moreover $(ii)$ of Assumption \ref{ass:reg:ext2} is satisfied by using that for all $p>s$
  \begin{equation*}\min_{1\leq m\leq\mdagn} \frac{\ev_{m}}{m(\hw_m)_{\vee1}}\asymp (\mdagn)^{-1-2a-(2s)_{\vee0}}  \asymp  n^{-(2a+1+ (2s)_{\vee0})/(2p+1 - 2s +(2a+2s)_{\vee0})} \geq 2 (\log n)/n .\end{equation*}
Finally, consider $(iii)$ of Assumption \ref{ass:reg:ext2}. It is easily verified that $ N_n \asymp n^{1/(1+(2s)_{\vee0})}$ which satisfies $\max_{1\leq m\leq N_n} \hw_{m} \leq N_n^{(2s)_{\vee0}}\asymp n^{(2s)_{\vee0}/(1+(2s)_{\vee0})} \leq  n$  and $M_n\asymp n^{1/(2a+1+(2s)_{\vee0})}\leq N_n \leq n$. Thereby also $(iii)$ of Assumption \ref{ass:reg:ext2} holds true.
\paragraph{Case [E-P].} We have  $M_n\asymp
n^{1/(2a+1+(2s)_{\vee 0})}$, $\mdagn\asymp \{\log [n (\log n)^{-(2a+1)/(2p)}] \}^{1/(2p)}$ and $ N_n \asymp n^{1/(1+(2s)_{\vee0})}$. Then as in case [P-P] $(i)$ and $(iii)$ of Assumption \ref{ass:reg:ext2} hold true since $M_n$ and $N_n$ are unchanged. Furthermore, for all $s\in\R$ we have
 \begin{equation*}\min_{1\leq m\leq\mdagn} \frac{\ev_{m}}{m(\hw_m)_{\vee1}}\asymp (\mdagn)^{-1-2a-(2s)_{\vee0}}  \asymp  (\log n)^{-(2a+1+ (2s)_{\vee0})/(2p)} \geq 2 (\log n)/n,\end{equation*} 
which shows (ii) of Assumption \ref{ass:reg:ext2}.
\paragraph{Case [P-E].} Here we have $M_n\asymp (\log  \frac{n\;(\log\log n)/(2a)}{(\log
n)^{(1+2a+(2s)_{\vee0})/(2a)}})^{1/(2a)} = (\log
n)^{1/(2a)}(1+o(1))$,  $\mdagn\asymp (\log \frac{n\;(\log\log n)/(2a)}{(\log n)^{(1+2a+2p)/(2a)}})^{1/(2a)} = (\log n)^{1/(2a)}(1+o(1))$ and $ N_n \asymp n^{1/(1+(2s)_{\vee0})}$. It is easily seen that  $(iii)$ of Assumption \ref{ass:reg:ext2} is satisfied.
Moreover, $(i)$ of Assumption \ref{ass:reg:ext2} holds true, since $\min_{1\leq j\leq  M_n} \ev_{j}\asymp \exp(-M_n^{2a})\asymp
\frac{(\log n)^{(1+2a+(2s)_{\vee0})/(2a)}}{n(\log\log n)/(2a)}\geq
2/n$ and 
\begin{equation*}\max_{m\geq  M_n} \frac{\ev_{m}}{m(\hw_m)_{\vee1}}\asymp
M_n^{-1-(2s)_{\vee0}}\exp(-M_n^{2a})\asymp \frac{(\log
n)}{n(\log\log n)/(2a)}\leq (\log n)/(2n).\end{equation*}
Finally, consider  $(ii)$ of Assumption \ref{ass:reg:ext2} which is satisfied by using that for all $p>s$
  \begin{equation*}\min_{1\leq m\leq\mdagn} \frac{\ev_{m}}{m(\hw_m)_{\vee1}}\asymp (\mdagn)^{-1-(2s)_{\vee0}}\exp(-(\mdagn)^{2a})  \asymp  \frac{(\log n)^{(2a+2p-(2s)_{\vee0})/(2a)}}{n\;(\log\log n)/(2a)}
     \geq 2 (\log n)/n,\end{equation*}
which completes the proof of Corollary \ref{optimadapt}.\hfill$\square$

\bibliography{FLM-SOP-MS}

\begin{thebibliography}{21}
\providecommand{\natexlab}[1]{#1}
\providecommand{\url}[1]{\texttt{#1}}
\expandafter\ifx\csname urlstyle\endcsname\relax
  \providecommand{\doi}[1]{doi: #1}\else
  \providecommand{\doi}{doi: \begingroup \urlstyle{rm}\Url}\fi

\bibitem[Barron et~al.(1999)Barron, Birg{\'e}, and Massart]{BBM99}
A.~Barron, L.~Birg{\'e}, and P.~Massart.
\newblock Risk bounds for model selection via penalization.
\newblock \emph{Probab. Theory Related Fields}, 113\penalty0 (3):\penalty0
  301--413, 1999.

\bibitem[Cardot and Johannes(2009)]{CardotJohannes2007}
H.~Cardot and J.~Johannes.
\newblock Thresholding projection estimators in functional linear models.
\newblock \emph{forthcoming in the Journal of Multivariate Analysis}, 2009.

\bibitem[Cardot et~al.(2003)Cardot, Ferraty, and Sarda]{CardotFerratySarda2003}
H.~Cardot, F.~Ferraty, and P.~Sarda.
\newblock Spline estimators for the functional linear model.
\newblock \emph{Statistica Sinica}, 13:\penalty0 571--591, 2003.

\bibitem[Comte et~al.(2006)Comte, Rozenholc, and Taupin]{CRT2006}
F.~Comte, Y.~Rozenholc, and M.-L. Taupin.
\newblock Penalized contrast estimator for density deconvolution.
\newblock \emph{Canadian Journal of Statistics}, 37\penalty0 (3), 2006.

\bibitem[Crambes et~al.(2009)Crambes, Kneip, and Sarda]{CrambesKneipSarda2007}
C.~Crambes, A.~Kneip, and P.~Sarda.
\newblock Smoothing splines estimators for functional linear regression.
\newblock \emph{Annals of Statistics}, 37\penalty0 (1):\penalty0 35--72, 2009.

\bibitem[Engl et~al.(2000)Engl, Hanke, and Neubauer]{EHN00}
H.~W. Engl, M.~Hanke, and A.~Neubauer.
\newblock \emph{Regularization of inverse problems.}
\newblock Kluwer Academic, Dordrecht, 2000.

\bibitem[Ferraty and Vieu(2006)]{FerratyVieu2006}
F.~Ferraty and P.~Vieu.
\newblock \emph{Nonparametric Functional Data Analysis: Methods, Theory,
  Applications and Implementations.}
\newblock Springer-Verlag, London, 2006.

\bibitem[Forni and Reichlin(1998)]{ForniReichlin1998}
M.~Forni and L.~Reichlin.
\newblock Let's get real: A factor analytical approach to disaggregated
  business cycle dynamics.
\newblock \emph{Review of Economic Studies}, 65:\penalty0 453--473, 1998.

\bibitem[Hall and Horowitz(2007)]{HallHorowitz2007}
P.~Hall and J.~L. Horowitz.
\newblock Methodology and convergence rates for functional linear regression.
\newblock \emph{Annals of Statistics}, 35\penalty0 (1):\penalty0 70--91, 2007.

\bibitem[James et~al.(2009)James, Wang, and Zhu]{AOStoappear}
G.~M. James, J.~Wang, and J.~Zhu.
\newblock Functional linear regression that's interpretable.
\newblock Technical report, To appear in the Annals of Statistics., 2009.

\bibitem[Johannes(2009)]{Johannes2009}
J.~Johannes.
\newblock Nonparametric estimation in circular functional linear model.
\newblock Technical report, University Heidelberg (revised and submitted),
  2009.
\newblock URL \url{http://arxiv.org/abs/0901.4266v1}.

\bibitem[Mair and Ruymgaart(1996)]{MairRuymgaart96}
B.~A. Mair and F.~H. Ruymgaart.
\newblock {Statistical inverse estimation in Hilbert scales.}
\newblock \emph{SIAM Journal on Applied Mathematics}, 56\penalty0 (5):\penalty0
  1424--1444, 1996.

\bibitem[Massart(2007)]{Massart07}
P.~Massart.
\newblock \emph{Concentration inequalities and model selection}, volume 1896 of
  \emph{Lecture Notes in Mathematics}.
\newblock Springer, Berlin, 2007.
\newblock Lectures from the 33rd Summer School on Probability Theory held in
  Saint-Flour, July 6--23, 2003, With a foreword by Jean Picard.

\bibitem[{M\"uller} and {Stadtm\"uller}(2005)]{MullerStadtmuller2005}
H.-G. {M\"uller} and U.~{Stadtm\"uller}.
\newblock Generalized functional linear models.
\newblock \emph{Ann. Stat.}, 33:\penalty0 774--805, 2005.

\bibitem[Natterer(1984)]{Natterer84}
F.~Natterer.
\newblock {Error bounds for Tikhonov regularization in Hilbert scales.}
\newblock \emph{Applicable Analysis}, 18:\penalty0 29--37, 1984.

\bibitem[Neumann(1997)]{Neumann1997}
M.~H. Neumann.
\newblock On the effect of estimating the error density in nonparametric
  deconvolution.
\newblock \emph{Journal of Nonparametric Statistics}, 7:\penalty0 307--330,
  1997.

\bibitem[Petrov(1995)]{Petrov1995}
V.~V. Petrov.
\newblock \emph{Limit theorems of probability theory. Sequences of independent
  random variables.}
\newblock Oxford Studies in Probability. Clarendon Press., Oxford, 4. edition,
  1995.

\bibitem[Preda and Saporta(2005)]{PredaSaporta2005}
C.~Preda and G.~Saporta.
\newblock Pls regression on a stochastic process.
\newblock \emph{{Computational Statistics \& Data Analysis}}, 48:\penalty0 149
  --158, 2005.

\bibitem[Ramsay and Silverman(2005)]{RamsaySilverman2005}
J.~Ramsay and B.~Silverman.
\newblock \emph{Functional Data Analysis.}
\newblock Springer, New York, second ed. edition, 2005.

\bibitem[Talagrand(1996)]{Tala96}
M.~Talagrand.
\newblock New concentration inequalities in product spaces.
\newblock \emph{Invent. Math.}, 126\penalty0 (3):\penalty0 505--563, 1996.

\bibitem[Tautenhahn(1996)]{Tautenhahn96}
U.~Tautenhahn.
\newblock {Error estimates for regularization methods in Hilbert scales.}
\newblock \emph{SIAM Journal on Numerical Analysis}, 33\penalty0 (6):\penalty0
  2120--2130, 1996.

\end{thebibliography}
\end{document}